\documentclass[12pt]{article}

\usepackage{amsthm,amsmath,amsfonts,epsfig,amssymb,latexsym, bbm}
\usepackage{mathrsfs} 
\usepackage[T1]{fontenc}
\usepackage{graphicx}
\usepackage{enumerate}
\usepackage{xy,xypic}
\usepackage{graphics}
\usepackage{footnote}
\usepackage{scalefnt}
\usepackage{tikz}
\usepackage{color}
\usetikzlibrary{shapes,arrows,fit,calc,positioning}
\usetikzlibrary{matrix}
\xyoption{all}
\usepackage{imakeidx}
\makeindex

\title{D\'{e}vissage Hermitian Theory} 
\author{Satya Mandal \\ University of Kansas,  Lawrence, Kansas 66045, USA }

\begin{document}

\pagenumbering{roman}
\setcounter{page}{0}

\renewcommand{\baselinestretch}{1.255}
\setlength{\parskip}{1ex plus0.5ex}
\date{15 August 2024} 
\newcommand{\iso}{\stackrel{\sim}{\longrightarrow}}
\newcommand{\sur}{\twoheadrightarrow}

\newcommand{\eop}{\hfill \TCP{\rule{2mm}{2mm}}}
\newcommand{\pf}{\noindent{\bf Proof.~}}
\newcommand{\outlinePf}{\noindent{\bf Outline of the Proof.~}}

\newcommand{\PD}{\dim_{{\SA}}}
\newcommand{\PDV}{\dim_{{\SV(X)}}}

\newcommand{\ra}{\rightarrow}
\newcommand{\lra}{\longrightarrow}
\newcommand{\hra}{\hookrightarrow} 
\newcommand{\Lra}{\Longrightarrow}
\newcommand{\Lla}{\Longleftarrow}
\newcommand{\Llra}{\Longleftrightarrow}
\newcommand{\llra}{\longleftrightarrow}
\newcommand{\pic}{The proof is complete.~}
\newcommand{\Sp}{\mathrm{Sp}}
\newcommand{\BiSp}{\mathrm{BiSp}}

\newcommand{\Dia}{\diagram}

\newcommand{\bE}{\begin{enumerate}}
\newcommand{\eE}{\end{enumerate}}

\newtheorem{theorem}{Theorem}[section]
\newtheorem{proposition}[theorem]{Proposition}
\newtheorem{lemma}[theorem]{Lemma}
\newtheorem{definition}[theorem]{Definition}
\newtheorem{corollary}[theorem]{Corollary}
\newtheorem{construction}[theorem]{Construction}
\newtheorem{notation}[theorem]{Notation}
\newtheorem{notations}[theorem]{Notations}
\newtheorem{remark}[theorem]{Remark}
\newtheorem{question}[theorem]{Question}
\newtheorem{example}[theorem]{Example} 
\newtheorem{examples}[theorem]{Examples} 
\newtheorem{exercise}[theorem]{Exercise} 
\newtheorem{clarification}[theorem]{Clarification}

\newtheorem{problem}[theorem]{Problem} 
\newtheorem{conjecture}[theorem]{Conjecture} 

\newcommand{\bD}{\begin{definition}}
\newcommand{\eD}{\end{definition}}
\newcommand{\bP}{\begin{proposition}}
\newcommand{\eP}{\end{proposition}}
\newcommand{\bL}{\begin{lemma}}
\newcommand{\eL}{\end{lemma}}
\newcommand{\bT}{\begin{theorem}}
\newcommand{\eT}{\end{theorem}}
\newcommand{\bC}{\begin{corollary}}
\newcommand{\eC}{\end{corollary}} 
\newcommand{\TCP}{\textcolor{purple}}
\newcommand{\TCM}{\textcolor{magenta}}
\newcommand{\TCR}{\textcolor{red}}
\newcommand{\TCB}{\textcolor{blue}}
\newcommand{\TCG}{\textcolor{green}}
\def\spec#1{\mathrm{Spec}\left(#1\right)}
\def\proj#1{\mathrm{Proj}\left(#1\right)}
\def\supp#1{\mathrm{Supp}\left(#1\right)}
\def\Sym#1{{\CS}\mathrm{ym}\left(#1\right)} 

 \def\LRf#1{\left( #1\right)}
\def\LRs#1{\left\{ #1\right\}}
\def\LRt#1{\left[ #1\right]}

\def\a{\mathfrak {a}}
\def\b{\mathfrak {b}}
\def\c{\mathfrak {c}}

\def\d{\mathfrak {d}}
\def\e{\mathfrak {e}}
\def\f{\mathfrak {f}}
\def\g{\mathfrak {g}}
\def\i{\mathfrak {i}}
\def\j{\mathfrak {j}}
\def\k{\mathfrak {k}}
\def\l{\mathfrak {l}}
\def\m{\mathfrak {m}}
\def\n{\mathfrak {n}}
\def\p{\mathfrak {p}}
\def\q{\mathfrak {q}}
\def\r{\mathfrak {r}}
\def\s{\mathfrak {s}}
\def\t{\mathfrak {t}}
\def\u{\mathfrak {v}}
\def\w{\mathfrak {w}}

\def\A{\mathfrak {A}}
\def\B{\mathfrak {B}}
\def\C{\mathfrak {C}}
\def\D{\mathfrak {D}}
\def\E{\mathfrak {E}}

\def\F{\mathfrak {F}}

\def\G{\mathfrak {G}}
\def\H{\mathfrak {H}}
\def\I{\mathfrak {I}}
\def\J{\mathfrak {J}}
\def\K{\mathfrak {K}}
\def\L{\mathfrak {L}}
\def\M{\mathfrak {M}}
\def\N{\mathfrak {N}}
\def\O{\mathfrak {O}}
\def\P{\mathfrak {P}}

\def\Q{\mathfrak {Q}}

\def\R{\mathfrak {R}}

\def\Sf{\mathfrak {S}} 
\def\T{\mathfrak {T}}
\def\U{\mathfrak {U}}
\def\V{\mathfrak {V}}
\def\W{\mathfrak {W}}

\def\CA{\mathcal {A}}
\def\CB{\mathcal {B}}
\def\CP{\mathcal {P}}
\def\CC{\mathcal {C}}
\def\CD{\mathcal {D}}
\def\CE{\mathcal {E}}
\def\CF{\mathcal {F}}
\def\CG{\mathcal {G}}
\def\CH{\mathcal{H}}
\def\CI{\mathcal{I}}
\def\CJ{\mathcal{J}}
\def\CK{\mathcal{K}}
\def\CL{\mathcal{L}}
\def\CM{\mathcal{M}}
\def\CN{\mathcal{N}}
\def\CO{\mathcal{O}}
\def\CP{\mathcal{P}}
\def\CQ{\mathcal{Q}}

\def\CR{\mathcal{R}}

\def\CS{\mathcal{S}}
\def\CT{\mathcal{T}}
\def\CU{\mathcal{U}}
\def\CV{\mathcal{V}}
\def\CW{\mathcal{W}} 
\def\CX{\mathcal{X}}
\def\CY{\mathcal{Y}}
\def\CZ{\mathcal{Z}}

\newcommand{\smallcirc}[1]{\scalebox{#1}{$\circ$}}

\def\BA{\mathbb{A}}
\def\BB{\mathbb{B}}
\def\BC{\mathbb{C}}
\def\BD{\mathbb{D}}
\def\BE{\mathbb{E}}
\def\BF{\mathbb{F}}
\def\BG{\mathbb{G}}
\def\BH{\mathbb{H}}
\def\BI{\mathbb{I}}
\def\BJ{\mathbb{J}}
\def\BK{\mathbb{K}}
\def\BL{\mathbb{L}}
\def\BM{\mathbb{M}}
\def\BN{\mathbb{N}} 
\def\BO{\mathbb{O}}
\def\BP{\mathbb{P}} 
\def\BQ{\mathbb{Q}}
\def\BR{\mathbb{R}}
\def\BS{\mathbb{S}}
\def\BT{\mathbb{T}}
\def\BU{\mathbb{U}}
\def\BV{\mathbb{V}}
\def\BW{\mathbb{W}}
\def\BX{\mathbb{X}}
\def\BY{\mathbb{Y}}
\def\BZ{\mathbb {Z}}

\def\SA{\mathscr {A}} 
\def\SB{\mathscr {B}}
\def\SC{\mathscr {C}}
\def\SD{\mathscr {D}}
\def\SE{\mathscr {E}}
\def\SF{\mathscr {F}}
\def\SG{\mathscr {G}}
\def\SH{\mathscr {H}}
\def\SI{\mathscr {I}}
\def\SK{\mathscr {K}} 
\def\SL{\mathscr {L}} 
\def\SM{\mathscr {M}}
\def\SN{\mathscr {N}}
\def\SO{\mathscr {O}}
\def\SP{\mathscr {P}}
\def\SQ{\mathscr {Q}}
\def\SR{\mathscr {R}}
\def\SS{\mathscr {S}}
\def\ST{\mathscr {T}}
\def\SU{\mathscr {U}}
\def\SV{\mathscr {V}}
\def\SW{\mathscr {W}}
\def\SX{\mathscr {X}}
\def\SY{\mathscr {Y}}
\def\SZ{\mathscr {Z}} 
 
\def\bfA{\bf A}
\def\bfB{\bf B}
\def\bfC{\bf C}
\def\bfD {\bf D}
\def\bfE{\bf E}
\def\bfF{\bf F}
\def\bfG{\bf G}
\def\bfH{\bf H}
\def\bfI{\bf I}
\def\bfJ{\bf J}
\def\bfK{\bf K}
\def\bfL{\bf L}
\def\bfM{\bf M}
\def\bfN{\bf N} 
\def\bfO{\bF O}
\def\bfP{\bF P} 
\def\bfQ{bBF Q}
\def\bfR{\bF R}
\def\bfS{\bF S}
\def\bfT{\bF T}
\def\bfU{\bf U}
\def\bfV{\bf V}
\def\bfW{\bf W}
\def\bfX{\bf X}
\def\bfY{\bf Y}
\def\bfZ{\bf Z}

\maketitle



\pagenumbering{arabic}

\noindent{\bf Abstract:} We prove D\'{e}vissage theorems for Hermitian $K$ Theory
(or $GW$ theory), analogous to Quillen's D\'{e}vissage theorem for $K$-theory. 
For abelian categories ${\SA}:=\left({\SA}, ^{\vee}, \varpi\right)$  with duality, and appropriate abelian subcategories ${\SB}\subseteq {\SA}$, we prove D\'{e}vissage theorems for 
${\bf GW}$ spaces, $G{\mathcal W}$-spectra and ${\mathbb G}W$ bispectra. As a consequence, for 
regular local rings $(R, \m, \kappa)$ with $1/2\in R$, we compute the ${\mathbb G}W$ groups 
${\mathbb G}W^{[n]}_k\left(\spec{R}\right)~\forall k, n\in {\mathbb Z}$, where $n$ represent the translation. 
\section{Introduction} 
Following \cite{QSS79}, for a category ${\SA}$, and object 
$M$ in ${\SA}$ will be referred to as ${\SA}$-module. Given an injective morphism $K \hra M$ (resp. {\it surjective morphism} $M\sur N$) , we say $K$ is a submodule  
(resp. quotient modules) of $M$. 

 In this article we  discuss 
D\'{e}vissage theorems for hermitian $K$-theory, also known as Grothendieck-Witt ($GW$) theory. For abelian  categories  D\'{e}vissage theorem for ${\bfK}$-theory spaces was proved by Quillen (\ref{introDeviss}), as follows. %

\bT[Quillen]\label{introDeviss}{\rm 
Let ${\SA}$ be a small abelian category, and ${\SB}\subseteq {\SA}$ be a full subcategory. Assume that ${\SB}$ is closed under submodules, quotient modules and finite products in ${\SA}$. 
Further, assume that every ${\SA}$ module $M$ has a finite filtration by ${\SB}$modules
(see Def \ref{DefFilll}).
Then the induced map
$$
{\bfK}\LRf{{\SB}} \lra {\bfK}\LRf{{\SA}}  
\quad {\rm is~ a~ homotopy~equivalence~of~spaces.}
$$
Equivalently, the maps $K_k\LRf{{\SB}} \iso K_k\LRf{{\SA}}$ of the $K$-groups are isomorphisms $\forall~k\geq 0$. 
}
\eT 
See \cite{Q}, \cite[pp 107]{M23} for more details.
Under the additional hypothesis that ${\SA}$ is noetherian,
the negative groups ${\BK}_{-k}\LRf{\SA}=0~\forall k\geq 1$ \cite[Thm 7]{S06}. 
 Consequently, the map of ${\BK}$-theory spectra
 $$
{\BK}\LRf{{\SB}} \lra {\BK}\LRf{{\SA}}  
\quad {\rm is~ a~ homotopy~equivalence,}
$$
in the category $\Sp$ of spectra of pointed topological spaces. 
By incorporating duality, in this article, we prove variety of D\'{e}vissage theorems 
of Hermitian theory, namely,  (1) D\'{e}vissage ${\bf GW}$ spaces, (2) D\'{e}vissage $G{\CW}$ spectra, and (3) D\'{e}vissage Karoubi ${\BG}{W}$  bispectra, for abelian categories ${\SA}=\LRf{{\SA}, ^{\vee}, \varpi}$ with duality and appropriate subcategories ${\SB}\subseteq {\SA}$, as in the setup (\ref{SetupDev}), similar to that in (\ref{introDeviss}).
For  abelian categories ${\SA}=\LRf{{\SA}, ^{\vee}, \varpi}$ with duality, the $GW$ spaces ${\bf GW}({\SA})$ and $GW$ spectrum $G{\CW}({\SA})\in \Sp$  are directly defined \cite{S10}, \cite[pp 424, 426]{M23}. However, $\forall n\in {\BZ}$, the Karoubi shifted ${\BG}W^{[n]}({\SA}) \in 
\BiSp$ bi spectrum is defined as 
$$
{\BG}W^{[n]}({\SA}):={\BG}W^{[n]}\LRf{{\bf dg}{\SA}, \q, ^{\vee}, \varpi}
\qquad \forall n\in {\BZ}. 
$$
where $\LRf{{\bf dg}{\SA}, \q, ^{\vee}, \varpi}$ denotes the dg category of ${\SA}$ with quasi isomorphism $\q$ and induced duality $^{\vee}, \varpi$ \cite[Sec 8]{S17}, \cite[pp 465]{M23}. The D\'{e}vissage Theorem for  the Karoubi ${\BG}W$ bispectrum  (\ref{BGWDeviss}) states that, when $1/2\in {\SA}$ and ${\SB}\subseteq {\SA}$ is a subcategory, as in the setup
 (\ref{SetupDev}), similar to that in (\ref{introDeviss}), the natural map
 $$
 {\BG}W^{[n]}\LRf{{\SB}, ^{\vee}, \varpi} \iso  {\BG}W^{[n]}\LRf{{\SA}, ^{\vee}, \varpi} 
 \qquad in \qquad \BiSp
 $$
 is a homotopy equivalence in $\BiSp, ~\forall n\in {\BZ}$. Consequently, we have natural isomorphisms
 $$
 {\BG}W^{[n]}_k\LRf{{\SB}} \iso  {\BG}W^{[n]}_k\LRf{{\SA}} \qquad \forall n, k\in {\BZ}
 $$
of the $ {\BG}W$ groups. 

The motivating drive for this article was the structure of the ${\BG}W\LRf{C{\BM}^Z(\spec{R})}$ where $R=(R, \m, \kappa)$ is a regular local ring, $1/2\in R$. 
Here $\LRf{C{\BM}^Z(\spec{R})}$ denotes the category of $R$modules $M$ of finite length (with $\PDV(M)=d:=\dim R$. 
 We prove 
(\ref{DivCMmXpcmi})
$$
  {\BG}{W}^{[n]}_k\LRf{C{\BM}^Z(X)} \cong 
  \left\{\begin{array}{ll}
  W\LRf{{\SV}(\kappa), ^{\vee}, \varpi} &  k\leq -1, n-k-d=0~mod~4\\ 
  0 &  k\leq -1, ~n-k-d=1, 2, 3~mod~4\\ 
  G{W}_k\LRf{{\SV}(\kappa), ^{\vee}, \varpi} & n=0~mod ~4,~ k\geq 0\\
  G{W}_k\LRf{{\SV}(\kappa), ^{\vee}, -\varpi} & n=2~mod ~4, k\geq 0\\
    \end{array}\right.
 $$
 This structure was instrumental in computing the ${\BG}W$-structure \cite{M24} of a 
 punctured subscheme $U=X-\{\m\}$ of a quasi projective scheme $X$, where $\m$ is a regular point. 
Here and throughout, $W(-)$ denotes  
Witt groups,  $GW(-)$ denotes Grothendieck-Witt groups.

{\it I am tankful to Amalendu Krishna 
for  useful communications and encouragements. }


\section{Preliminaries} %
For variety of standard definitions we refer to \cite{M23, QSS79, Q, S10}. We will mainly be interest in abelian  categories with duality, while most of the definitions make sense for exact categories with duality. Recall the definitions of
Exact category $\LRf{{\SE}, ^{\vee}, \varpi}$ \cite[pp 421]{M23}. We quickly recall some of the basic definitions and notations.

\bD\label{DefBasic}{\rm 
Suppose $\LRf{{\SE}, ^{\vee}, \varpi}$ is an exact category with duality. We alway assume $\varpi$ is an equivalence and $2$ is invertible in ${\SE}$. 
An object $M\in Obj\LRf{{\SE}}$ may be referred to as and {\bf ${\SE}$-module}. 
For an injective map $\i: L\subseteq M$, we may say $L$ is a {\bf submodule} or a {\bf subobject} of $M$. 
A {\bf symmetric/hermitian form} is a pair $(M, \varphi)$ where $M\in Obj({\SE})$ and 
$\varphi:M\lra M^{\vee}$ is a morphism, such that $\varphi=\varpi_M \varphi^{\vee}$. 
A {\bf symmetric/hermitian space} is a symmetric form $(M, \varphi)$, such that $\varphi$ is an isomorphism.  Given a  hermitian space $(M, \varphi)$ and submodule 
$\i: L\subseteq M$, we write $L^{\perp}=\ker\LRf{\i^{\vee}\varphi}$, to be called 
the {\bf orthogonal complement} of $L$. Note that there is no inclusion or disjoiness relationship between $L$ and $L^{\perp}$. However, if $L\subseteq L^{\perp}$, 
equivalently if
$\iota^{\vee}\varphi \iota=0$,  then we say that $L$ is {\bf totally isotropic} submodule of $M$. 

For any small category ${\SC}$, its classifying space is denoted by ${\BB}{\SC}$ \cite[pp 36]{M23}.
 For a small exact category ${\SE}$ the ${\BQ}$-category of ${\SE}$ will be denoted by 
${\BQ}{\SE}$ \cite[pp 61]{M23}. So, the ${\bfK}$-theory space of ${\SE}$ is defined as
${\bfK}\LRf{{\SE}}=\Omega {\BB}{\BQ}{\SE}$, the loop space \cite[pp 124]{M23}. 
}
\eD 

We improvise the  the setup of Quillen \cite[Def 3.6.1 pp 107]{M23}, for D\'{e}vissage theorem for $K$-theory spaces (\ref{introDeviss}), to suit out the context of abelian  categories with duality.

\bD[Setup]\label{SetupDev}{\rm 
Let $\LRf{{\SA}, ^{\vee}, \varpi}$ be an abelian category with duality. This means 
\bE
\item First, ${\SA}$ is an abelian category.
\item There is an exact functor (duality)
$
^{\vee}: {\SA}^{op} \lra {\SA}
$.
\item There is a double dual identification 
$
\varpi:1_{\SA} \lra ^{\vee\vee}
$, which is a natural equivalence. 
\item Further, in this context of Hermitian theory, we assume $1/2\in {\SA}$. This means, for ${\SA}$-modules $M, N$, the map $2:Hom(M, N) \sur Hom(M, N)$ is surjective.
\eE 
\noindent{\bf Setup:}
Let ${\SB}\subseteq {\SA}$ be a full  subcategory, as in 
\cite[Def 3.6.1 pp 107]{M23}. So:\\
(1) ${\SB}$ is closed under submodules, quotient modules and finite product.\\
(2) Consequently, ${\SB}$ is closed under finite sum, in a module in ${\SA}$. If $B_1, \ldots, B_k\in {\SB}$, and are submodules of $M\in {\SA}$ then $B_1+B_2+\cdots +B_k\in {\SB}$. \\
Further, assume that ${\SB}$  closed under duality. This means,
\bE
\item $B\in Obj({\SB}) \Lra B^{\vee}\in Obj({\SB})$.
\item $\forall B\in Obj({\SB})$ the isomorphism $\varpi_B$ is in ${\SB}$, which is true because ${\SB}$ is a full subcategory. 

\eE 
}
\eD
We will be working under this setup, for the rest of this article. 
\subsection{Jordan-H\"{o}lder Theorem} 
\vspace{3mm}
We recall the following  definitions of filtrations  
 from \cite[pp 285]{QSS79}. 
 
 \bD\label{DefFilll}{\rm
 Consider the setup (\ref{SetupDev}). 
 \bE
 \item\label{602OneFIL} Let $N$ be an ${\SA}$-module. A ${\SB}$-filtration of $N$ is a finite sequence of submodules
 $$
 \diagram
 0=N_0 \ar@{^(->}[r] & N_1  \ar@{^(->}[r]  & \cdots 
  \ar@{^(->}[r] & N_{r-1}  \ar@{^(->}[r] & N_r=N\\
 \enddiagram
 ~ \ni ~ \frac{N_i}{N_{i-1}} \in {\SB} ~\forall i=1, 2, \ldots, r
 $$
 \item Let $(M, \varphi)$ be a hermitian space in ${\SA}$.  A hermitian ${\SB}$-filtration consists of a totally isotropic submodule $N\subseteq M$  together with a ${\SB}$-filtration of $N$, as above (\ref{602OneFIL}),
  such that $\frac{N^{\perp}}{N} \in {\SB}$. 
  Consequently, there is ${\SB}$-filtration of $M$ as follows:
  $$
 \diagram
 0=N_0 \ar@{^(->}[r] & N_1  \ar@{^(->}[r]  & \cdots 
  \ar@{^(->}[r] & N_{r-1}  \ar@{^(->}[r] & N_r \ar@{^(->}[r] &\\
  N_r^{\perp}  \ar@{^(->}[r] &N_{r-1}^{\perp}\ar@{^(->}[r] &\cdots \ar@{^(->}[r] &N_1^{\perp} \ar@{^(->}[r] &N_0^{\perp}&M\ar@{=}[l]\\
 \enddiagram
 $$
 \item It is worth pointing out  to the readers that any two hermitian ${\SB}$-filtrations  of $(M, \varphi)$ have a common refinement \cite[6.6 Thm]{QSS79}. 
 \eE 
} \eD
 
 For completeness, we include the following defintion.
 \bD{\rm 
 Let ${\SA}$ be an abelian category and $M$ be an ${\SA}$-module. We say that $M$ is noetherian if every increasing sequence of submodules $M$ terminates. We say that 
 ${\SA}$ is noetherian if every ${\SA}$-modules $M$ is noetherian. 
 }
 \eD

 \vspace{3mm} 
\bL\label{socleFor}{\rm 
Consider the setup (\ref{SetupDev}). Also assume that every $M\in {\SA}$ has a finite 
${\SB}$-filtration. 
Further, assume that ${\SA}$ is noetherain.
%
Then 
\bE
\item\label{507MaxBmod} Given $0\neq M\in {\SA}$, there is a unique maximal ${\SB}$-submodule, to be denoted by $\m(M)$, such that  
$0\neq \m(M)\subseteq M$. If $M=0$ write $\m(M)=0$. 

 \item Let $(M, \varphi)$ be a  Hermitian space in ${\SA}$. Denote 
$$
\s(M)=\m(M) \cap \m(M)^{\perp}\subseteq M
$$
Then 
$$
\s(M)=0 \Llra M\in {\SB}
$$
\eE
}
\eL
\pf Given an ${\SA}$-module $M\neq 0$, it has nonzero submodules, because of the filtration condition. Therefore,  $M$ has maximal ${\SB}$-submodules because of the noetherian 
condition. It is unique,  because ${\SB}$ is closed under sum. This establishes the first proposition. 

To prove the later proposition,
suppose $M\in {\SB}$.
Then $\m(M)=M$ and hence $\m(M)^{\perp}=0$. So, $\s(M)=0$. Conversely,
suppose $\s(M)=0$. Write $N=\m\LRf{\m(M)^{\perp}}\in {\SB}$. So, 
$N\subseteq \m(M)$ and $N\subseteq \m(M)^{\perp}$. So, $N\subseteq \s(M)=\m(M) \cap \m(M)^{\perp}=0$. By (\ref{507MaxBmod}), we have $\m(M)^{\perp}=0$. So,
$M=\m(M)\in {\SB}$.  
\pic $\eop$ 

\vspace{3mm} 
The following is a Jordan-H\"{o}lder Theorem, which is in the spirit of \cite{QSS79}.

\bL[Jordan-H\"{o}lder Theorem]\label{Jordan}{\rm 
Consider the setup of (\ref{SetupDev}).  Assume that %
 every ${\SA}$-module $N$ %
 has a finite  ${\SB}$-filtration:
\begin{equation}\label{UsuusQFi}
0=N_0\subseteq N_1\subseteq \cdots \subseteq N_r=N
\quad \ni \left\{\begin{array}{l}
N_i\in Obj\LRf{{\SA}} ~\forall i=0, 1, \ldots, r\\
\frac{N_i}{N_{i-1}}\in Obj\LRf{{\SB}} ~\forall i= 1, \ldots, r\\
\end{array}\right.
\end{equation}
Further assume that ${\SA}$ is noetherain.
%
%
Let $(M, \varphi)$ be a hermitian space in ${\SA}$. Then $(M, \varphi)$ has a {\bf hermitian ${\SB}$-filtration}. Due to hypothesis (\ref{UsuusQFi}), this only requires that $M$ has a totally isotropic submodule $N\subseteq M$ such that 
$\frac{N^{\perp}}{N}\in {\SB}$. Consequently, we would have the following:

\bE
\item Given such a totally isotopic submodule $N=:N_r$ 
with $\frac{N^{\perp}}{N}\in {\SB}$, by hypothesis \ref{UsuusQFi}, there is a finite ${\SB}$-filtration:
$$
0=N_0 \subseteq N_1 \subseteq N_2 \subseteq \cdots \subseteq N_r=N \qquad 
\ni \quad \frac{N_i}{N_{i-1}}\in {\SB}\quad \forall i= 1, \ldots, r
$$

\item Further, $\frac{N_{i-1}^{\perp}}{N_i^{\perp}}\cong \LRf{\frac{N_i}{N_{i-1}}}^{\vee}
\in {\SB}$. %
\item \label{MainFilt14} Since, $\frac{N_r^{\perp}}{N_r}\in {\SB}$, 
\begin{equation}\label{JurEfnfHH}
0=N_0 \subseteq N_1 \subseteq N_2 \subseteq \cdots \subseteq N_{r-1}
\subseteq N_r
\subseteq  N_r^{\perp} \subseteq N_{r-1}^{\perp} \subseteq \cdots 
\subseteq N_1^{\perp} \subseteq N_0^{\perp}=M
\end{equation}
is a  ${\SB}$-filtration of $M$. 
\eE 
}
\eL
\pf We only need to prove the existence of totally isotropic submodule $N$, with $\frac{N^{\perp}}{N}\in {\SB}$. We assume $M\neq 0$. Let $N_0=0$.
\bE
\item By Lemma \ref{socleFor}, $M$ has  a maximal ${\SB}$-submodule $0\neq \m(M) \subseteq M$. 
Write  $L_1:=\m(M)\in {\SB}$. 
Write $N_1=\s(L_1)=L_1\cap L_1^{\perp}$, which is totally isotropic. 
Further $N_1 \subseteq L_1^{\perp} \subseteq N_1^{\perp}$. 
%
%
 Note $N_1\in {\SB}$, because ${\SB}$ is  closed under subobjects. 
By %
\cite[5.2]{QSS79}, $\varphi$ induces a natural nonsingular form:
$$
\overline{\varphi}_1: \frac{N_1^{\perp}}{N_1} \lra \LRf{\frac{N_1^{\perp}}{N_1}}^{\vee}
$$
\item Now let $0\neq L_2=\m\LRf{\frac{N_1^{\perp}}{N_1}}\subseteq  \frac{N_1^{\perp}}{N_1}$ be the maximal submodule in ${\SB}$. Let $\s(L_2):=L_2\cap L_2^{\perp}\subseteq \frac{N_1^{\perp}}{N_1}$. Let $N_1\subseteq N_2\subseteq N_1^{\perp}$ be such that 
$\frac{N_2}{N_1}=\s(L_2) \subseteq \frac{N_1^{\perp}}{N_1}$. Note 
$\frac{N_2}{N_1}\in {\SB}$. In fact $N_2$ is totally isotropic in $M$, by \cite[6.5 Prop (2)]{QSS79}.

\bE
\item 
Suppose $\frac{N_2}{N_1}=\s(L_2)=0$. In this case $N_1=N_2$. By 
Lemma \ref{socleFor}, $\frac{N_1^{\perp}}{N_1}\in {\SB}$. 

\item Assume $\frac{N_2}{N_1}=\s(L_2)\neq 0$. In this case $N_1\neq N_2$.
By Lemma \ref{socleFor} $\frac{N_1^{\perp}}{N_1}\notin {\SB}$. Further, $\frac{N_2}{N_1}$ is totally isotropic in $\LRf{\frac{N_1^{\perp}}{N_1}, \overline{\varphi}_1}$. 
Again, by %
\cite[5.2]{QSS79}
$\overline{\varphi}_1$ induces a natural nonsingular form 
$$
\diagram
 \frac{\s(L_2)^{\perp}}{\s(L_2)}\ar[rr]^{\overline{\varphi}_2}\ar[d]_{\wr} &&\LRf{\frac{\s(L_2)^{\perp}}{\s(L_2)}}^{\vee}\\
\frac{N_2^{\perp}}{N_2} \ar[rr]&&\LRf{\frac{N_2^{\perp}}{N_2}}^{\vee}\ar[u]^{\wr}\\
\enddiagram
$$
\eE 
\item The process terminates because of the noetherian condition.
So, we obtain a sequence of totally isotropic submodules 
$$
\diagram 
N_0=0 \ar@{^(->}[r] & N_1 \ar@{^(->}[r] &N_2 \ar@{^(->}[r] &\cdots \ar@{^(->}[r] &N_{r-1} \ar@{^(->}[r] &N_r \ar@{^(->}[r] & \cdots  \ar@{^(->}[r] &M\\
\enddiagram 
$$
Such that 
$$
L_{r+1}=\m\LRf{\frac{N_{r}^{\perp}}{N_{r}}}\subseteq \frac{N_{r}^{\perp}}{N_{r}}, \quad \s(L_{r+1})=L_{r+1}\cap L_{r+1}^{\perp}=0
$$
By Lemma \ref{socleFor} $\frac{N_{r}^{\perp}}{N_{r}}\in {\SB}$. Further, 
$\frac{N_{r+1}}{N_{r}}=\s(L_r)=0$ for some $N_{r+1}\subseteq M$. 
Hence $N_{r+1}=N_{r}$

\eE 
Since $\frac{N_{r}^{\perp}}{N_{r}}\in {\SB}$.
\pic $\eop$ 

\vspace{3mm}
As an immediate consequence of Jordan-H\"{o}rder theorem \ref{Jordan}, we retrieve the 
D\'{e}vissage theorem for Witt groups, under the setup \ref{SetupDev}. For an abelian category $\LRf{{\SA}, ^{\vee}, \varpi}$ with duality, we refer to \cite{QSS79}, \cite[pp 423]{M23} for the definition of the Witt group $W\LRf{{\SA}, ^{\vee}, \varpi}$.

\bC[Witt D\'{e}vissage]\label{WittDev618}{\rm 
Consider the setup (\ref{SetupDev}).  Assume ${\SA}$ is noetherian. Then the following maps 
$$
\left\{\begin{array}{l}
W\LRf{{\SB}, ^{\vee}, \varpi}\iso W\LRf{{\SA}, ^{\vee}, \varpi}\\
W\LRf{{\SB}, ^{\vee}, -\varpi}\iso W\LRf{{\SA}, ^{\vee}, -\varpi}\\
\end{array}\right.
\quad {\rm are~isomorphisms.} 
$$
}
\eC
\pf This was proved in \cite[6.9 Cor]{QSS79}, under   \cite[Hypothesis {\bf $B_0$}]{QSS79}, which is same as validity of  Lemma \ref{Jordan}. \pic $\eop$ 

\vspace{3mm} 
For future reference, we recall some of the results on Balmer Witt groups
\cite{B00, B01, BW02}.
\begin{remark}\label{notaBalTriW}{\rm 
For a  Triangulated category  $\LRf{{\SK}, \#, T, \varpi}$ with duality, with $1/2\in {\SK}$, Balmer Witt groups $W\LRf{{\SK}, \#, T, \varpi}$ was defined in \cite{B00}, \cite[pp 503]{M23}. Further, for 
$n\in {\BZ}$ the $n$-shifted Witt group was defined as 
$W^n\LRf{{\SK}, \#, T, \varpi}:=W\LRf{T^n\LRf{{\SK}, \#, T, \varpi}}$.
Now,  let $\LRf{{\SA}, ^{\vee}, \varpi}$, be an abelian category, with duality and $1/2\in {\SA}$.
 \bE
 \item Define
 $$
 W^n\LRf{{\SA}, ^{\vee}, \varpi}:=W^n\LRf{{\bf D}^b\LRf{{\SA}}}\qquad \qquad \forall n\in {\BZ}
 $$ 
 These groups are referred to as Balmer Witt groups. 
 \item We have
 $$
  W^n\LRf{{\SA}, ^{\vee}, \varpi}=
 \left\{\begin{array}{lll}
0 & if~n=2r+1~odd& by ~\cite[Prop ~5.2]{BW02}\\
W\LRf{{\SA}, ^{\vee}, \varpi} & if~n=0 ~mod~4n&by ~ \cite[4.7 ~Thm ~pp~ 376]{B01}\\
W\LRf{{\SA}, ^{\vee}, -\varpi} & if~n=2 ~mod~4n&by~\cite[pp ~378]{B01}\\
\end{array}\right.
  $$
  \item Using this equation, it follows from (\ref{WittDev618}), under the hypotheses of (\ref{WittDev618}), the natural maps
  $$
W^n\LRf{{\SB}, ^{\vee}, \varpi}  \iso W^n\LRf{{\SA}, ^{\vee}, \varpi} 
\quad {\rm are~isomorphisms}\qquad \forall n\in {\BZ}. 
  $$
  \item\label{pointFour} Let $\LRf{R, \m, \kappa}$ be a regular local commutative ring, with $d=\dim R$. Write $X=\spec{R}$, $Z=V(\m)$. Let $C{\BM}^Z(X)$ denote the category of $R$-modules of finite length.  
   For 
  $M\in C{\BM}^Z(X)$, the association $M^{\vee}\mapsto Ext^d\LRf{M, R}$ defines a natural duality on  $C{\BM}^Z(X)$. It is a Theorem in \cite[Thm 6.1]{BW02} that 
  $$
  W^n\LRf{C{\BM}^Z(X), ^{\vee}, \varpi}=
  \left\{\begin{array}{ll}
  W(\kappa) & n=\dim R ~mod~4\\
  0 & otherwise\\
  \end{array}\right. 
  $$
 \eE 
}
\end{remark}

\section{D\'{e}vissage ${\bf GW}$ spaces and $G{\CW}$ spectra} 
First, we establish some preliminaries. 
We  recall the construction of the hermitian ${\BQ}$-category  ${\BQ}^h\LRf{{\SE}, ^{\vee}, \varpi}$, of an exact category $\LRf{{\SE}, ^{\vee}, \varpi}$ with duality, which is due to Giffen and Karoubi \cite[pp 12]{S10}, \cite[pp 424 Def 10.1.8]{M23}.
\bD\label{inQHMasd}{\rm 
Let ${\SE}=\LRf{{\SE}, ^{\vee}, \varpi}$ be a small exact category with duality. 
Let ${\BQ}{\SE}$ denote the ${\BQ}$-category of ${\SE}$ \cite[pp 61]{M23}.
Define the 
category ${\BQ}^h\LRf{{\SE}, ^{\vee}, \varpi}$ as follows:
\bE
\item Objects of ${\BQ}^h\LRf{{\SE}, ^{\vee}, \varpi}$ are the symmetric spaces $(P, \varphi)$.
\item A morphism $(P_1, \varphi_1) \lra (P_2, \varphi_2)$ is a morphism $\psi:=
(E, p, \iota): P_1\lra P_2$ in the category ${\BQ}{\SE}$, 
such that the diagram
\begin{equation}\label{392DefQHerObj}
\diagram 
P_1 \ar[d]_{\varphi_1}^{\wr}&& E\ar@{-->}[d]^{\chi} \ar@{->>}[ll]_p \ar@{^(->}[rr]^{\iota} & &P_2\ar[d]^{\varphi_2}_{\wr}\\
P_1^{\vee} \ar@{^(->}[rr]_{p^{\vee}}&& E^{\vee}   && P_2^{\vee}\ar@{->>}[ll]^{\iota^{\vee}}\\
\enddiagram
\quad \left\{\begin{array}{l}
commutes.~ It ~follows  \\
\ker(\iota^{\vee}\varphi)=\LRf{\ker(p), \iota \zeta}\\
\end{array}\right.
\end{equation}
\eE
}
\eD
 So, by definition 
$$
Mor_{{\BQ}\LRf{{\SE}, ^{\vee}, \varpi}}\LRf{(P_1, \varphi_1), (P_2, \varphi_2)} \subseteq 
Mor_{{\BQ}{\SE}}\LRf{P_1, P_2} 
$$

\vspace{3mm}
We define the hermitian  analogue of {\bf admissible layers}, \cite[pp 68]{M23} due to Quillen.
\bD\label{AmdLyQP}{\rm 
Let $\LRf{{\SE}, ^{\vee}, \varpi}$ be an exact category with duality. Let $(P, \varphi)$ be a quadratic space. A totally isotropic submodule 
$N \subseteq P$ is defined as an {\bf admissible layer} in $(P, \varphi)$.
Let 
\begin{equation}\label{ISOSusse}
{\SI}so(P, \varphi)=\LRs{N\subseteq P: {\rm where ~N~ is ~a~ totally ~isotropic~ submodule ~of }~ (P, \varphi)}
\end{equation} 
Note ${\SI}so(P, \varphi)$ is a partially ordered set, by superset. So, $N_1\leq N_2
\Llra N_2 \subseteq N_1$. So, ${\SI}so(P, \varphi)$  is treated as a  category.

In the linear case, an admissible layer in a module $M$ is defined as a pair of submodules $K_0\subseteq K_2 \subseteq M$ \cite[pp 68]{M23}.
Let 
$$
{\CS}(M)=\LRs{
(K_0, K_1): K_0\subseteq K_1\subseteq M ~{\rm is~an~admissible~layer~in~}M
}
$$
This is also a partially ordered set. 
There  is natural order preserving map
$$
{\SI}so(P, \varphi) \lra {\CS}(P)\quad \qquad {\rm sending}\quad N \mapsto \LRf{N, N^{\perp}}
$$
}
\eD 
The following is an equivalent way to characterize the maps in 
$(P_1, \varphi_1) \lra (P_2, \varphi_2)$ in ${\BQ}^h\LRf{{\SE}, ^{\vee}, \varpi}$, with admissible layers in $(P_2, \varphi_2)$, as clarified below. 
\bL\label{One0ftu414}{\rm 
The maps (\ref{392DefQHerObj}) are in bijection, with the set of admissible layers
$K\in {\SI}so(P, \varphi)$, as follows 
\begin{equation}\label{410EquIValal}
\diagram
K \ar@/^/@{^(->}[rd] \ar@{^(->}[d] & \\
E \ar@{^(->}[r] & P_2\\
\enddiagram
\ni \left\{\begin{array}{ll}
K\subseteq K^{\perp} \subseteq P_2& {\rm i.e~K~is~totally~isotropic}\\
K^{\perp}=E& \frac{E}{K}\cong P_1\\
By~general~theory ~\varphi_2~induces& \overline{\varphi}:\frac{E}{K}\iso \LRf{\frac{E}{K}}^{\vee}\\
 \LRf{\frac{E}{K}, \overline{\varphi}} \cong \LRf{P_1, \varphi} & is~an~isometry\\
\end{array}\right.
\end{equation}
}
\eL
\pf Consider a map  (\ref{392DefQHerObj}). Write $K:=\ker(p)$. 
We have $\iota^{\vee}\varphi_2\iota\zeta=0$. So, $(\iota\zeta)^{\vee}\varphi_2\iota\zeta=0$. So, $\ker(p) \subseteq \ker(p)^{\perp}$, or 
$K \subseteq K^{\perp}$. It is obvious that $\frac{E}{K} \cong P_1$.
 Further, 
$$
\zeta^{\vee}\chi= \zeta^{\vee}(p^{\vee}\varphi_1p)
=  (p\zeta)^{\vee}\varphi_1p=0.
\quad So, \quad E \subseteq K^{\perp}.
$$
Denote the inclusion maps $\eta:E \hra K^{\perp}$, $\beta: K^{\perp}\hra P_2$. Consider the diagram
$$
\diagram 
P_1\ar[rrrr]^{\varphi_1}&&&&P_1^{\vee}\ar@{^(->}[dd]^{p^{\vee}}&&\\
~~K^{\perp} 
\ar@/^/[rrrru]^{\gamma}\ar@/^/@{_(->}[rrd]_{\beta} \ar@/^/[rrrrd]^{\iota^{\vee}\varphi _2\beta}&&&&&&\\
\ar@/^/@{->>}[uu]^{p}E\ar@/_/@{^(->}[u]_{\eta}  \ar@{^(->}[rr]_{\iota}  && P_2\ar@{->>}[rr]_{\iota^{\vee}\varphi_2} && E^{\vee} \ar@{->>}[rr]_{\zeta^{\vee}}&& K^{\vee} \\
\enddiagram 
Note\quad \left\{\begin{array}{l}
\LRf{P_1^{\vee}, p^{\vee}}=\ker\LRf{\zeta^{\vee}}\\
\zeta^{\vee}\iota^{\vee}\varphi_2\beta=0\\
So, ~\exists~\gamma~ \ni ~\p^{\vee}\gamma=\iota^{\vee}\varphi_2\beta\\
\end{array}\right.
$$
We have
$$
(p^{\vee}\gamma)\eta =\iota^{\vee} \varphi_2 \beta\eta =\iota^{\vee} \varphi _2\iota
= p^{\vee} \psi p\Lra p^{\vee}\left(\gamma \eta- \varphi_1p\right)=0
$$
So, $\gamma \eta=\varphi_1p$. Therefore, $\gamma \eta$ is surjective. So, $\gamma$ is surjective. Further, 
$$
\ker(\gamma)=\ker(\iota^{\vee}\varphi_2\beta)
=\ker(\iota^{\vee}\varphi_2)\cap K^{\perp}=K\cap K^{\perp}= K
$$
Going modulo $K$ we have the commutative diagram: 
$$
\diagram 
P_1 \ar[rr]_{\sim}^{\varphi_1} && P_1^{\vee}\\
\frac{E}{K} \ar[rr]_{\overline{\eta}} \ar[u]_{\wr}&&\frac{K^{\perp}}{\ker(p)}\ar[u]_{\wr}\\
\enddiagram
$$
So, $\overline{\eta}$ is an isomorphism, since other three maps are isomorphisms. It follows from this that $\eta$ is isomorphism (use Snake). Therefore, it is established that $E=K^{\perp}$.
By general theory, $\varphi_2$ induces hermitian space
$$
\overline{\varphi}: \frac{E}{K} \iso \LRf{\frac{E}{K}}^{\vee}
$$
It is easy to check
$$
\diagram
 P_1\ar[d]_{\varphi_1} & \frac{E}{K} \ar[d]^{\overline{\varphi}}\ar[l]_{\sim}\\
P_1^{\vee}\ar[r]_{\sim} &\LRf{\frac{E}{K}}^{\vee}\\
 \enddiagram
 \quad {\rm is~an~isometry.} 
$$
So, the Data is (\ref{410EquIValal}) is established. 

Conversely, assume we have a set of data as in (\ref{410EquIValal}). 
Existence of $\overline{\varphi}$ follows from general theory. Since 
$\overline{p}:\LRf{\frac{E}{K}}\iso P_1$ is an isomorphism, it induces 
a hermitian space $P_1\iso P_1^{\vee}$ structure on $P_1$. So,  the diagram

$$ %
\diagram 
P_1 \ar[d]_{\varphi_1}^{\wr}&& E\ar@{-->}[d]^{\chi} \ar@{->>}[ll]_p \ar@{^(->}[rr]^{\iota} & &P_2\ar[d]^{\varphi_2}_{\wr}\\
P_1^{\vee} \ar@{^(->}[rr]_{p^{\vee}}&& E^{\vee}   && P_2^{\vee}\ar@{->>}[ll]^{\iota^{\vee}}\\
\enddiagram
\qquad {\rm commutes.}
$$
So, we obtained a map $(P_1, \varphi_1) \lra (P_2, \varphi_2)$ in ${\BQ}^h({\SE})$.  
\pic 
$\eop$

\vspace{3mm} 
The following corollary provides a good handle for our subsequent discussions. 
\bC\label{MapstoADM}{\rm 
Let ${\SE}=\LRf{{\SE}, ^{\vee}, \varpi}$ be a small exact category with duality. 
Let $(P, \varphi)$ be a Hermitian space. 
As usual, let $1_{{\BQ}^h\LRf{{\SE}, ^{\vee}, \varpi}}/(P, \varphi)$ 
denote the category of all maps $(L, \psi) \lra (P, \varphi)$ in 
${\BQ}^h\LRf{{\SE}, ^{\vee}, \varpi}$, where $(L, \psi)$ runs through all the hermitian spaces. 
Let ${\SI}so(P, \varphi)$ be as defined in (\ref{ISOSusse}), as a category.
Then there is an equivalence of categories 
$$
{\SI}so(P, \varphi) \lra 1_{{\BQ}^h\LRf{{\SE}, ^{\vee}, \varpi}}/(P, \varphi) 
\quad sending\quad N \mapsto \LRf{
\LRf{\frac{N^{\perp}}{N}, \varphi_N}\lra (P, \varphi)}
$$
where $\varphi_N:\frac{N^{\perp}}{N} \iso \LRf{\frac{N^{\perp}}{N}}^{\vee}$ is the induced hermitian space. 
}
\eC
\pf Follows from Lemma \ref{One0ftu414}. $\eop$

 
\subsection{D\'{e}vissage ${\BQ}^h\LRf{{\SA, ^{\vee}, \varpi}}$} \label{SECDev2QQheGW}

In this section we first prove an analog to Quillen's D\'{e}vissage \cite[pp 107 Def 3.6.1]{M23}, for the category ${\BQ}^{h}\LRf{{\SE}, ^{\vee}, \varpi} $. Subsequently, we use this to prove the analogue of the same for Grothendieck-Witt (${\bf GW}$) spaces.
%

\vspace{3mm}
\bL{\rm 
Consider the setup of (\ref{SetupDev}). 
Then 
\bE
\item\label{OneQDR} The natural functor
$$
{\BQ}{\SB} \lra {\BQ}{\SA} \quad {\rm is~a~full~subcategory}. 
$$
\item\label{ZwiQTY} Further, the natural functor 
$$
{\BQ}^h\LRf{{\SB}, ^{\vee}, \varpi} \lra {\BQ}^h\LRf{{\SA}, ^{\vee}, \varpi} \quad {\rm is~a~full~subcategory}. 
$$
\eE 
}
\eL
\pf The first assertion was proved in the paper of  Quillen \cite{Q}, \cite[Def 3.6.1]{M23}. So, we prove the second one.
We have a commutative diagram
$$
\diagram
Mor_{{\BQ}^h\LRf{{\SB}, ^{\vee}, \varpi}}\LRf{(P_1,\varphi_1), (P_2, \varphi_2)} \ar@{^(->}[d]\ar[rr] && 
Mor_{{\BQ}^h\LRf{{\SA}, ^{\vee}, \varpi}}\LRf{(P_1, \varphi_1), (P_2, \varphi_2)}\ar@{^(->}[d]
 \\
Mor_{{\BQ}{\SB}}\LRf{P_1, P_2} \ar[rr]_1 && Mor_{{\BQ}{\SA}}\LRf{P_1, P_2} \\
\enddiagram 
$$
Since ${\BQ}{\SB}$ is a full subcategory of  ${\BQ}{\SA}$  by (\ref{OneQDR}), 
 the lower arrow is a bijection. So, the upper arrow in injective. 
To see fullness of the upper arrow, let $(p, U, \iota): (P_1, \varphi_1) \lra (P_2, \varphi_2)$ be a map in ${\BQ}^h\LRf{{\SA}, ^{\vee}, \varpi}$, while $(P_1, \varphi_1)$, $(P_2, \varphi_2)$ are in $\LRf{{\SB}, ^{\vee}, \varpi}$. So, we have the diagram
$$
\diagram 
&&\ker(p) \ar@{^(->}[d]\ar[rr]^{\iota_0}_{\sim}&&\ker(\iota^*\varphi_2)\ar@{^(->}[d]\\
P_1 \ar[d]_{\varphi_1}^{\wr}&& U\ar@{-->}[d] \ar@{->>}[ll]_p \ar@{^(->}[rr]^{\iota} & &P_2\ar[d]^{\varphi_2}_{\wr}\\
P_1^{\vee} \ar@{^(->}[rr]_{p^*}&& U^{\vee}   && P_2^*\ar@{->>}[ll]^{\iota^{\vee}}\\
\enddiagram
\quad where\quad \iota_0:=\iota_{|\ker(p^{\vee})}
$$
The object $U$ and the maps  $\iota$ and $p$ are in ${\SA}$.
However, since ${\SB}$ is closed under submodules, $U\in Obj({\SB})$. Similarly, 
$ker(p), \ker(\iota^*\varphi_2)\in Obj({\SB})$. Now, by fullness of ${\SB}$ in ${\SA}$ all maps are in ${\SB}$. So, the map $(p, U, \iota)$ is in ${\BQ}\LRf{{\SB}, ^{\vee}, \varpi}$.
\pic $\eop$ 


\vspace{3mm}
We insert the following lemma, before we proceed.
\bL\label{conttble}{\rm 
Consider the setup of (\ref{SetupDev}), and the inclusion functor
$$
\f:{\BQ}\LRf{{\SB}, ^{\vee}, \varpi}\hra {\BQ}\LRf{{\SA}, ^{\vee}, \varpi}
$$
Let $\LRf{M, \varphi}\in {\SB}$.
Then $\f/\LRf{M, \varphi}$ is contractible.
}
\eL 
\pf As usual, $\LRf{M, \varphi}$ has an identity element in ${\BQ}^h\LRf{{\SA}, ^{\vee}, \varpi}$, namely
$$
1_{\LRf{M, \varphi}}:=\qquad \qquad 
\diagram
M\ar[d]_{\varphi}  & M\ar[d]^{\varphi} \ar@{^(->}[r]^1 \ar@{->>}[l]_1 & M\ar[d]^{\varphi}\\
M^{\vee} & M^{\vee}\ar@{^(->}[l]^1 \ar@{->>}[r]_1  & M^{\vee}\\
\enddiagram
$$
Since $(M, \varphi)\in {\BQ}^h\LRf{{\SB}, ^{\vee}, \varpi}$, we have 
$$
\f/\LRf{M, \varphi} = 1_{ \LRf{{\BQ}^h\LRf{{\SB}, ^{\vee}, \varpi}}}/\LRf{M, \varphi} 
$$
which has a terminal element, namely $1_{\LRf{M, \varphi}}$.
So, $\f/{\LRf{M, \varphi}}$ is contractible, by usual  generalities \cite[pp 40]{M23}.
\pic $\eop$  


\vspace{3mm}
The following is 
a  theorem of Schlichting \cite[Thm 6.1]{S10}, while we substitute the hypothesis (b) by the noetherian
condition (Theorem \ref{Jordan}). To conform with our own style, and for the benefit of the readers, we include the proof. 
\bT[Schlichting]\label{740MainDididi}{\rm 
Consider the setup of (\ref{SetupDev}).  
Assume that every ${\SA}$-module $M$ has a finite ${\SB}$-filtration.
Further, assume ${\SA}$ is noetherian.
%
Then the inclusion functor 
\begin{equation}\label{incluDev}
\f:{\BQ}\LRf{{\SB}, ^{\vee}, \varpi}\hra {\BQ}\LRf{{\SA}, ^{\vee}, \varpi}
\end{equation}
 induces a homotopy equivalence
 $$
{\BB}{\BQ}\LRf{{\SB}, ^{\vee}, \varpi}\lra {\BB}{\BQ}\LRf{{\SA}, ^{\vee}, \varpi}
$$
of the classifying spaces. 
}
\eT

\pf Consider the inclusion functor $\f$, as in (\ref{incluDev}). 
{\bf Fix any  
$(M, \varphi) \in Obj{\BQ}^h\LRf{{\SA}, ^{\vee}, \varpi}$}. So, $(M, \varphi)$ is a symmetric space in 
$\LRf{{\SA}, ^{\vee}, \varpi}$. 
By Theorem A
of Quillen \cite[pp 49]{M23}, it is enough to prove that 
 the category $\f/(M, \varphi)$ is contractible. By (\ref{Jordan}) there is a
 hermitian ${\SB}$-filtration
\begin{equation}\label{NumSimbly}
0=N_0 \subseteq N_1 \subseteq N_2 \subseteq \cdots \subseteq N_{r-1}
\subseteq N_r
\subseteq  N_r^{\perp} \subseteq N_{r-1}^{\perp} \subseteq \cdots 
\subseteq N_1^{\perp} \subseteq N_0^{\perp}=M
\end{equation}
where $\forall i~N_i$ are totally isotropic in $(M, \varphi)$ and   $\frac{N_r^{\perp}}{N_r}\in {\SB}$. Let $\varphi_i:\frac{N_i^{\perp}}{N_i}\iso \LRf{\frac{N_i^{\perp}}{N_i}}^{\vee}$ denote the symmetric space induced by $\varphi$ (which is same as the one induced by $\varphi_{i-1}$). 
%
%
The category $\f/\LRf{\frac{N_r^{\perp}}{N_r}, \varphi_r}$ is contractible, by 
 Lemma \ref{conttble}.

{\bf Fix $0\leq k \leq r-1$.} 
 We write $M_k=\frac{N_{k}^{\perp}}{N_{k}}$ and $\LRf{M_k, \varphi_k}$ the induced hermitian space.
 Given a hermitian ${\SB}$-filtration, as above (\ref{NumSimbly}), there is a natural map
$$
\Phi_k: \LRf{M_{k+1}, \varphi_{k+1}} \lra \LRf{M_{k}, \varphi_{k}} 
\qquad {\rm in }\qquad \LRf{{\BQ}^h\LRf{{\SA}, ^{\vee}, \varpi}}
$$
as follows
$$
\Phi_k:=
\diagram
&\frac{N_{k+1}}{N_{k}}\ar@{^(->}[d]\ar[r]^1&\frac{N_{k+1}}{N_{k}}\ar@{^(->}[d]\\
 \frac{N_{k+1}^{\perp}}{N_{k+1}}\ar[d]_{\varphi_{k+1}} & \frac{N_{k+1}^{\perp}}{N_{k}} \ar[d]\ar@{^(->}[r]^{\iota_k}\ar@{->>}[l]_{p_k} & 
 \frac{N_{k}^{\perp}}{N_{k}} \ar[d]^{\varphi_{k}}\\
  \LRf{\frac{N_{k+1}^{\perp}}{N_{k+1}}}^{\vee}\ar@{^(->}[r]_{p_k^{\vee}}  & \LRf{\frac{N_k^{\perp}}{N_{k+1}}}^{\vee}& \LRf{\frac{N_{k}^{\perp}}{N_{k}}}^{\vee} \ar@{->>}[l]^{\iota_k^{\vee}} \\
\enddiagram 
$$
 By composition, with $\Phi_k$, we obtain a functor
$$
F_k: \f/\LRf{M_{k+1}, \varphi_{k+1}} \lra \f/\LRf{M_{k}, \varphi_{k}} 
~ \diagram 
(L, \psi) \ar[r]^{f\quad} \ar@/_/[dr]& \LRf{ \frac{N_{k+1}^{\perp}}{N_{k+1}}, \varphi_{k+1}} \ar[d]^{\Phi_ko-}\\
& \LRf{ \frac{N_{k+1}^{\perp}}{N_{k}}, \varphi_{k}} \\
\enddiagram 
~ \left\{\begin{array}{l}
(L, \psi)\in {\BQ}^h\LRf{{\SB}, ^{\vee}, \varpi}\\
F_k(f)=\Phi_kof\\
\end{array}\right.
$$
  As pointed out above,  $\f/\LRf{M_r, \varphi_r}$ is contractible.  
  We will prove that $F_k$ is a homotopy equivalence $\forall k=0, 1, \ldots, k-1$; 
  and this will complete the proof that  $\f/\LRf{M, \varphi}$ is contractible. 
We prove only the  case $k=1$, because there is no loss of generality, if  $\LRf{M_k, \varphi_k}$ is repalced by
$\LRf{M_0, \varphi_0}=\LRf{M, \varphi}$.

We will write $N:=N_1$, $p:=p_1$, $\i:=\iota_1$, $\Phi=\Phi_1$.  So, we have a map
$$
\Phi:=
\diagram
&N\ar@{^(->}[d]\ar[r]^1&N\ar@{^(->}[d]\\
M_1\ar[d]_{\varphi_{1}} &N^{\perp} \ar[d]^{\chi}\ar@{^(->}[r]^{\i}\ar@{->>}[l]_{p} & 
M \ar[d]^{\varphi}\\
 M_1^{\vee}\ar@{^(->}[r]_{p^{\vee}}  & \LRf{N^{\perp}}^{\vee}& M^{\vee} \ar@{->>}[l]^{\i^{\vee}} \\
\enddiagram 
\quad {\rm with}\quad 
N\in {\SB}
$$
The functor $F_1$ is given by 
$$
F_1: \f/\LRf{M_{1}, \varphi_{1}} \lra \f/\LRf{M, \varphi} 
~~ \diagram 
(L, \psi) \ar[r]^f \ar@/_/[dr]& \LRf{M_1, \varphi_1}  
\ar[d]^{\Phi o-}\\
& \LRf{M, \varphi} \\ 
\enddiagram 
\quad \left\{\begin{array}{l}
 (L, \psi)\in {\BQ}^h\LRf{{\SB}, ^{\vee}, \varpi}\\
 \Phi(f)=\Phi o f\\
 \end{array}\right.
$$
Refer to Corollary \ref{MapstoADM}, and define a subset  
subset ${\SI}so\LRf{M, \varphi; {\SB}}\subseteq {\SI}so\LRf{M, \varphi}$,
as follows:
$$
{\SI}so\LRf{M, \varphi; {\SB}}=\LRs{
K\subseteq M: {\rm  ~K~ is ~a~ totally ~isotropic~ submodule ~of }~ (M, \varphi), 
\ni \frac{K^{\perp}}{K}\in {\SB}
}
$$
The set   ${\SI}so\LRf{M, \varphi; {\SB}}$ inherits the partial order, by superset. 
Consider the commutative diagram
$$
\diagram 
{\SI}so\LRf{M, \varphi; {\SB}} \ar[rr]^{\sim}\ar@{^(->}[d]  && \f/(M, \varphi)\ar@{^(->}[d]\\
{\SI}so\LRf{M, \varphi}\ar[rr]_{\sim} && 1_{{\BQ}^h\LRf{{\SA}, ^{\vee}, \varpi}}/(M, \varphi)\\
\enddiagram
$$
By Corollary \ref{MapstoADM}, the second horizontal functor is  an equivalence. 
It is easy to check that the top horizontal functor is also an equivalence. 
The following diagram
$$
\diagram 
{\SI}so\LRf{M_1, \varphi_1; {\SB}} \ar[rr]^{G} \ar[d]_{\wr}&& {\SI}so\LRf{M, \varphi; {\SB}}\ar[d]^{\wr}&K \mapsto p^{-1}K\supseteq N\\
\f/\LRf{M_1, \varphi_1}\ar[rr]_{F_1} && \f/\LRf{M, \varphi}&\star \mapsto \star o \Phi\\
\enddiagram
$$
commutes. Since the vertical functors are equivalence, we prove that $G$ is a homotopy equivalence. 
Define a functor 
$$
\diagram 
{\SI}so\LRf{M, \varphi; {\SB}} \ar[r]^{\rho}& {\SI}so\LRf{M_1, \varphi_1; {\SB}}\\
\enddiagram 
\quad \rho(K)=
\left\{\begin{array}{lll}
{\bf 0}=\frac{N}{N}\subseteq M_1& if ~N\not\subseteq K&K \not\leq N\\
\frac{K }{N} \subseteq M_1& if ~N\subseteq K&K \leq N\\
\end{array}\right.
$$ 
Note $N$ is totally isotropic. 
Therefore, if $K$ is totally isotropic in $M$, so is $\rho(K)$ in $M_1$. By \cite[6.5 Prop]{QSS79} $\frac{\rho(K)^{\perp}}{\rho(K)}=\frac{K^{\perp}}{K} \in {\SB}$. So, $\rho(K) \in {\SI}so\LRf{M_1, \varphi_1; {\SB}}$. 
Suppose 
$K_1 \leq K_2$.
\bE
\item In case $K_2\leq N$ then $N\subseteq K_2\subseteq K_1$, Hence 
$\frac{K_2}{N} \subseteq \frac{K_1}{N}$. So, $\rho(K_1) \leq \rho(K_2)$.
\item Suppose $K_2\not\leq N$. So, $N \not\subseteq K_2$, and $\rho(K_2)=\frac{N}{N}=0$. If $N\subseteq K_1$ then $\rho(K_1)=\frac{K_1}{N} \leq \rho(K_2)=\frac{N}{N}=0$.
If $N\not\subseteq K_1$ then $\rho(K_1)=0$ and $\rho(K_1) \leq \rho(K_2)$.
\eE
So, this association $K \mapsto \rho(K)$ respects the partial order. Therefore $\rho$ is a functor. It follows 
$$
\rho G\LRf{\frac{K}{N}}=\frac{K}{N}, \quad {\rm which ~means } \quad 
\rho G=1_{{\SI}so\LRf{M_1, \varphi_1; {\SB}}}
$$
We define another functor
$$
\diagram 
{\SI}so\LRf{M, \varphi; {\SB}} \ar[r]^{\s}& {\SI}so\LRf{M, \varphi; {\SB}}\\
\enddiagram \qquad  {\rm by}\qquad \s(K)=K\cap N
$$
Since $K$ is totally isotropic, so is $\s(K)$.  Consider the commutative diagram
$$
\diagram
K^{\perp}\bigoplus N^{\perp} \ar@{->>}[rr]  \ar@{->>}[d] && K^{\perp} + N^{\perp}\ar@{->>}[d] \\
\frac{K^{\perp}}{K}\bigoplus \frac{N^{\perp}}{N}   \ar@{-->>}[rr] && \frac{K^{\perp} + N^{\perp}}{K\cap N}\\
\enddiagram
$$
It follows, that the broken arrow is surjective. 

Since $\frac{K^{\perp}}{K}\bigoplus \frac{N^{\perp}}{N} \in {\SB}$, it follows 
$\frac{\s(K)^{\perp}}{\s(K)}\cong \frac{K^{\perp} + N^{\perp}}{K\cap N} \in {\SB}$.
Clearly, $\s$ respects the partial order and hence $\s$ is a functor. 
  We define two natural transformation:
\bE
\item There is a natural transformation 
$$
\U: 1_{{\SI}so\LRf{M, \varphi; {\SB}}}\lra \s \qquad  \qquad \U_K:=\LRf{ K \leq  K \cap N}
$$
\item We have 
$$
G\rho(K)= \left\{\begin{array}{lll} 
K &if~ K \leq N & N \subseteq K\\
N & otherwise & \\
\end{array}\right. 
\qquad {\rm and}\quad \s(K) =K\cap N
$$
So, $G\rho(K) \leq \s(K)$. This defines a natural transformation
$$
\V: G\rho \lra \s \qquad \qquad \V_K:=\LRf{G\rho(K) \leq \s(K)}
$$
We need to check the naturallity. Let $K_1\leq K_2$ in ${\SI}so\LRf{M, \varphi; {\SB}}$. 
\bE
\item Suppose $ K_2 \leq N$. This means $N\subseteq K_2 \subseteq K_1$. In this case
$$
G\rho(K_1) =K_1 \leq G\rho(K_2) =K_2. \quad So, \quad 
\diagram 
K_1\ar[d]_{G\rho(K_1\leq K_2)} \ar[r]^{\V_{K_1}}& K_1\cap N \ar[d]^{\s(K_1\leq K_2)} \\
K_2  \ar[r]_{\V_{K_2}}& K_2\cap N\\
\enddiagram 
$$
commutes. 
\item Suppose $K_2\not \leq N$ and $K_1 \leq N$. So, $N\subseteq K_1$ and 
$N \not\subseteq K_2$. So, 
$$
G\rho(K_1)=K_1, \quad G\rho(K_2)= N. \quad So, \quad 
\diagram 
K_1\ar[d]_{G\rho(K_1\leq K_2)} \ar[r]^{\V_{K_1}}& K_1\cap N \ar[d]^{\s(K_1\leq K_2)} \\
N  \ar[r]_{\V_{K_2}}& K_2\cap N\\
\enddiagram 
$$ 
commutes. 
\item Suppose $K_2\not \leq N$ and $K_1 \not\leq N$. So, $N\not\subseteq K_1$ and 
$N \not\subseteq K_2$. So, 
$$
G\rho(K_1)=N, \quad G\rho(K_2)= N. \quad So, \quad 
\diagram 
N\ar[d]_{G\rho(K_1\leq K_2)} \ar[r]^{\V_{K_1}}& K_1\cap N \ar[d]^{\s(K_1\leq K_2)} \\
N  \ar[r]_{\V_{K_2}}& K_2\cap N\\
\enddiagram 
$$ 
commutes. 

\eE 

\eE 
In conclusion, $\U$ gives  a homotopy equivalence from $1_{{\SI}so\LRf{}M, \varphi; {\SB}}$ to 
$\s$, and $\V$ gives a homotopy equivalence from $G\rho$ to $\s$. Therefore, $G\rho$ is equivalent to  $1_{{\SI}so\LRf{}M, \varphi; {\SB}}$. 
Since $\rho G =1$, it follows $G$ is a homotopy equivalence. 
This completes the proof of the theorem. $\eop$ 
\section{D\'{e}vissage ${\bf GW}$ spaces and $G{\CW}$ spectra}\label{SECDuoSpect}
In this subsection we establish the Hermitian analogues of Quillen's D\'{e}vissage theorem  \cite[pp 107]{M23}.
We recall the definition of the Grothedieck-Witt spaces ${\bf GW}\LRf{{\SE}, ^{\vee}, \varpi}$ 
from \cite[pp 14]{S10}, \cite[pp 424]{M23} 
\bD\label{DEFGw}{\rm 
Let ${\SE}=\LRf{{\SE}, ^{\vee}, \varpi}$ be a small exact category with duality. 
Define the Grothedieck-Witt space ${\bf GW}\LRf{{\SE}, ^{\vee}, \varpi}$,  by the following homotopy fibration:
$$
\diagram 
{\bf GW}\LRf{{\SE}, ^{\vee}, \varpi} \ar[r]&  {\BB}{\BQ}^h\LRf{{\SE}, ^{\vee}, \varpi} \ar[r] & 
 {\BB}{\BQ}\LRf{{\SE}}\\
\enddiagram
$$
 Here the second arrow is the natural continuous map induced by the forgetful functor 
 ${\BQ}^h\LRf{{\SE}, ^{\vee}, \varpi} \lra {\BQ}\LRf{{\SE}}$.
}
\eD

\vspace{3mm} 
The D\'{e}vissage theorem for ${\bf GW}$ spaces is an immediate consequence of the definition \ref{DEFGw}, as follows.
\bT[D\'{e}vissage ${\bf GW}$ space]\label{HermDev}{\rm 
Consider the setup of (\ref{SetupDev}). So, we have an abelian category
$\LRf{{\SA}, ^{\vee}, \varpi}$ with duality, and a subcategory $\LRf{{\SB}, ^{\vee}, \varpi}$,
as in (\ref{SetupDev}). 
Assume ${\SA}$ is noetherian, $1/2\in {\SA}$, and that every ${\SA}$-module $M$ has a finite 
${\SB}$-filtration. Then there is a natural homotopy equivalence 
$$
{\bf GW}\LRf{{\SB}, ^{\vee}, \varpi} \lra {\bf GW}\LRf{{\SA}, ^{\vee}, \varpi}.
$$
of the ${\bf GW}$-spaces. 
In particular, there are natural isomorphisms 
$$
{\bf GW}_n\LRf{{\SB}, ^{\vee}, \varpi} \cong {\bf GW}_n\LRf{{\SA}, ^{\vee}, \varpi} ~\forall n\geq 0\qquad  {\rm of}~{\bf GW}{\rm -groups.}
$$  
}
\eT
\pf By definition (\ref{DEFGw}) we have a commutative diagram of homotopy fibrations:
$$
\diagram 
{\bf GW}({\SB}) \ar[r]\ar@{-->}[d]&  {\BB}{\BQ}^h\LRf{{\SB}, ^{\vee}, \varpi} \ar[r] \ar[d]& {\BB}{\BQ}\LRf{{\SB}}\ar[d]\\
{\bf GW}({\SA}) \ar[r]&  {\BB}{\BQ}^h\LRf{{\SA}, ^{\vee}, \varpi} \ar[r] & {\BB}{\BQ}\LRf{{\SA}}\\
\enddiagram
$$
The two horizontal arrows on the right side are the natural maps,
induced by forgetful functors, of the classifying spaces. 
From the commutative diagram of the long exact sequences of the homotopy groups we obtain isomorphisms 
${\bf GW}_n({\SB}) \iso {\bf GW}_n({\SA})~\forall n\geq 0$. 
Therefore, by Whitehead's theorem \cite[pp 638]{M23}, the map
${\bf GW}({\SB}) \iso {\bf GW}({\SA})$ is a homotopy equivalence. 
\pic $\eop$ 

 \vspace{3mm}
 Now we prove the D\'{e}vissage theorem for $G{\CW}$-spectra. For an exact category 
 ${\SE}:=\LRf{{\SE}, ^{\vee}, \varpi}$ there are two definition of $G{\CW}$-spectra possible. These take values in the category $\Sp$ of spectra of pointed spaces. 
 \bE
 \item First, 
 for such an exact category, a spectrum  ${\SE} \mapsto G{\CW}({\SE}):=G{\CW}\LRf{{\SE}, ^{\vee}, \varpi}: \in \Sp$ was  associated, directly,  in \cite{S10}. This was done by defining suspension functor $S:{\bf eCat}D \lra {\bf eCat}D$, where ${\bf eCat}D$ denotes the category of small exact categories with duality.  Note, given $\LRf{{\SE}, ^{\vee}, \varpi}\in {\bf eCat}D$ we obtain two spctra
 $$
 G{\CW}\LRf{{\SE}, ^{\vee}, \varpi}, ~G{\CW}\LRf{{\SE}, ^{\vee}, -\varpi}\in \Sp.
 $$
 \item Further, 
 for a
 pretriangualted  dg category ${\SA}:=\LRf{{\SA}, \q,  ^{\vee}, \varpi}$ with weak equivalence and duality a 
 $G{W}$-spectra $G{\CW}({\SA})\in {\Sp}$ 
 was defined in \cite[pp 1776]{S17},  \cite[pp 378, 449]{M23}. 
\item  Our interest to dg categories, in this article, remains limited to 
 ${\bf dg}{\SA}$ of abelian categories ${\SA}$. Therefore, all the dg categories we discuss are pretraingulated. 
 Subsequently, for an exact category ${\SE}$ (considered as a ${\BZ}$-linear category) 
denote the corresponding dg category, of bounded chain complexes of ${\SE}$-modules,  by ${\bf dg}{\SE}$. Note that the underlying exact category of ${\bf dg}{\SE}$ is the category $Ch^b\LRf{{\SE}}$ of bounded chain complexes, and the corresponding triangulated category $\T\LRf{{\SE}}={\bf D}^b\LRf{{\SE}}$ is the derived category of ${\SE}$.  The quasi isomorphisms in 
${\bf dg}{\SE}$ will be denoted by $\q$. 
Assume $1/2\in {\SE}\in {\bf eCat}D$. So, we obtain a spoectrum
$$
G{\CW}\LRf{{\bf dg}{\SE}, \q, ^{\vee}, \varpi}\in {\Sp}, \quad {\rm with~four~periodicity.} 
$$
Incorporating the $n$-translations, we obtain spectra
$$
G{\CW}^{[n]}\LRf{{\bf dg}{\SE}, \q, ^{\vee}, \varpi} \in \Sp\qquad \forall n\in {\BZ}. 
$$

\item 
 In fact 
$$
\left\{\begin{array}{ll}
G{\CW}^{[n]}\LRf{{\bf dg}{\SE}, \q, ^{\vee}, \varpi}=
G{\CW}^{[n+4]}\LRf{{\bf dg}{\SE}, \q, ^{\vee}, \varpi} & \forall n\in {\BZ}\\
G{\CW}\LRf{{\SE}, ^{\vee}, \varpi}=G{\CW}^{[4n]}\LRf{{\bf dg}{\SE}, \q, ^{\vee}, \varpi} 
& \forall n\in {\BZ}\\
G{\CW}\LRf{{\SE}, ^{\vee}, -\varpi}=G{\CW}^{[4n+2]}\LRf{{\bf dg}{\SE}, \q, ^{\vee}, \varpi} & \forall n\in {\BZ}\\
\end{array}\right.
$$
\item As usual, denote the $G{\CW}$-groups
$$
\left\{\begin{array}{ll}
G{\CW}_k({\SE})=\pi_k\LRf{G{\CW}({\SE})} & \forall k\in {\BZ}\\
G{\CW}_k^{[n]}({\SE})=\pi_k\LRf{G{\CW}^{[n]}\LRf{{\bf dg}{\SE}, \q, ^{\vee}, \varpi}} & \forall k, n\in {\BZ}\\
\end{array}\right.
$$
\eE 
 
 \vspace{3mm}
 Now we D\'{e}vissage $G{\CW}$  spectra of exact categories with duality.
 \bT[D\'{e}vissage $G{\CW}$ Spectra]\label{HermSpDev}{\rm 
Consider the setup of (\ref{SetupDev}).
Assume ${\SA}$ is noetherian and that every ${\SA}$-module $M$ has a finite 
${\SB}$-filtration. Then the natural map 
$$
G{\CW}\LRf{{\SB}, ^{\vee}, \varpi} \lra G{\CW}\LRf{{\SA}, ^{\vee}, \varpi} \qquad in \qquad \Sp
$$
of spectra, is a homotopy equivalence in the category of spectra ${\Sp}$.
This means,   that the natural maps 
$$
G{\CW}_n\LRf{{\SB}, ^{\vee}, \varpi} \cong G{\CW}_n\LRf{{\SA}, ^{\vee},  \varpi} ~\forall n\in {\BZ}\qquad  {\rm of}~G{\CW}{\rm -groups~are~isomorphisms.}
$$ 

}
\eT
\pf We have \cite[pp 1781]{S17}, \cite[pp 450]{M23}
$$
G{\CW}\LRf{{\SA}, ^{\vee}, \varpi} \cong 
G{\CW}_n\LRf{{\bf dg}{\SA}, \q, ^{\vee}, \varpi} =
\left\{\begin{array}{ll}
{\bf GW}_n\LRf{{\SA}, ^{\vee}, \varpi} & n\geq 0\\
W^{-n}\LRf{{\SA}, ^{\vee}, \varpi} & n<0\\
\end{array}\right.
$$
and likewise for $G{\CW}_n\LRf{{\SB}, ^{\vee}, \varpi}$. 
It follows from Theorem \ref{HermSpDev} that 
${\bf GW}_n\LRf{{\SB}, ^{\vee}, \varpi} \cong {\bf GW}_n\LRf{{\SA}, ^{\vee}, \varpi} ~\forall n\geq 0$. Further, for $n<0$, from ( \ref{notaBalTriW}), we have
$$
W^{-n}\LRf{{\SA}, ^{\vee}, \varpi}=\left\{\begin{array}{lll}
0 & if~n=2r+1~odd& by ~\cite[Prop ~5.2]{BW02}\\
W\LRf{{\SA}, ^{\vee}, \varpi} & if~n=0 ~mod~4n&by ~ \cite[4.7 ~Thm ~pp~ 376]{B01}\\
W\LRf{{\SA}, ^{\vee}, -\varpi} & if~n=2 ~mod~4n&by~\cite[pp ~378]{B01}\\
\end{array}\right.
$$
Same formula applies for $\LRf{{\SB}, ^{\vee}, \varpi}$. Now, 
$$
\left\{\begin{array}{l}
W\LRf{{\SB}, ^{\vee}, \varpi} \cong W\LRf{{\SA}, ^{\vee}, \varpi}\\
W\LRf{{\SB}, ^{\vee}, -\varpi} \cong W\LRf{{\SA}, ^{\vee}, -\varpi}\\
\end{array}\right.
$$
by Witt D\'{e}vissage Corollary \ref{WittDev618}.
 \pic $\eop$

\section{D\'{e}vissage  Karoubi ${\BG}W$ bispectra} 
For a pretriangualted dg category $\LRf{{\SA}, \q, ^{\vee},  \varpi}$ with weak equivalence and duality, the Karoubi ${\BG}W$ bi Spectra, usually denoted by ${\BG}W\LRf{{\SA}}:={\BG}W\LRf{{\SA}, \q, ^{\vee},  \varpi} \in {\BiSp}$ takes 
value in the category $\BiSp:=\BiSp\LRf{{\BS}^1, \widetilde{{\BS}}^1_{\Sp}}$ of bi spectra, which is the category of symmetric   $\widetilde{{\BS}}^1_{\Sp}$ spectra
(see \cite[\S 8.2]{S17}, \cite[pp 465]{M23} for details). 
For a dg category $\LRf{{\SA},  \q, ^{\vee},\varpi}$ with weak equivalence $\q$ and duality, we denote 
\cite[pp 466]{M23} 
 $$
 \left\{\begin{array}{l}
  {\BG}{W}^{[n]}\LRf{{\SA}, ,\q, ^{\vee}  \varpi}\in {\BiSp}={\rm the}~{\BG}{W}~{\rm bi~spectra~of}~ {\SA}^{[n]}~~\forall ~{\rm shift}~n\in {\BZ} \\
  {\BG}{W}^{[n]}_k\LRf{{\SA}, ,\q,, ^{\vee}, \varpi} =\pi_k {\BG}{W}^{[n]}\LRf{{\SA}, ,\q,^{\vee},  \varpi}\in {\BiSp}
  ~~ \forall n, k\in {\BZ}~{\rm the ~homotopy~groups} \\
  {\rm Note,}~ {\BG}{W}^{[n]}\LRf{{\SA}, ,\q, ^{\vee}  \varpi} \cong 
   {\BG}{W}^{[n+4]}\LRf{{\SA}, ,\q, ^{\vee}  \varpi},~\forall n\in {\BZ}~{\rm has~periodicity~four}.\\
  \end{array}\right.
  $$
  
\vspace{3mm}

First, we compute some of some of the negative ${\BG}W$ groups. We will restrict ourselves to dg categories, which are pretriangualted. 
 
\bP\label{BGWandWitt335}{\rm 
Suppose ${\SA}$ is a {\bf pretriangualted dg category} \cite[pp 9]{S17}, \cite[pp 437]{M23} with duality and weak equivalence. Assume  that the 
underlying exact category $Z_0{\SA}$ is abelian and noetherian. Assume $1/2\in {\SA}$. Then
$$
{\BG}{W}_{-k}^{[n]}\LRf{{\SA}} \cong W^{n+k}\LRf{{\SA}} \qquad \forall n, k\in {\BZ},~ k\geq 1 
$$
Consequently,
$$
{\BG}{W}_{-(k+4)}^{[n]}\LRf{{\SA}} \cong W^{n+k}\LRf{{\SA}} \cong {\BG}{W}_{-k}^{[n]}\LRf{{\SA}} \cong {G}{\CW}_{-k}^{[n]}\LRf{{\SA}}  \qquad \forall n, k\in {\BZ},~k\geq 1 
$$
}
\eP 
\pf The last isomorphism in the second set of isomorphisms is a restatement of  \cite[Prop 6.3 pp 1785]{S17}. So, we prove only the first isomorphism. 
Fix $n\in {\BZ}$. 
Recall the exact triangle \cite[Thm 8.11 pp 1801]{S17} 
\begin{equation}\label{extDelts}
\diagram 
{\BG}W^{[n]}\LRf{\SA} \ar[r] & {\BK}\LRf{\SA} \ar[r] & {\BG}W^{[n+1]}\LRf{\SA} \ar[r] & 
S^1\wedge {\BG}W^{[n]}\LRf{\SA}\\
\enddiagram
\end{equation} 
in the triangulated category  $\BiSp$ of bi spectra. 
Before we write down the exact sequences of homotopy groups corresponding to this exact triangle, we point out two facts.
Since  $Z_0{\SA}$ is abelian, the associated triangulated category ${\T}{\SA}$ is idempotent complete 
\cite[2.8 Thm]{BS01}. Therefore, 
we have 
${\BG}W_k^{[n]}\LRf{{\SA}} ={G}W_k^{[n]}\LRf{{\SA}} ~\forall ~k\geq 0$ \cite[pp 1799]{S17}. 
Further, it is assumed that $Z_0{\SA}$ is  noetherian. It follows from  \cite[Thm 7]{S06} that ${\BK}_{-k}\LRf{{\SA}}=0~\forall k\geq 1$.  
%
So, $\forall n\in {\BZ}$, the exact triangle induces exact sequences
$$
\diagram 
\cdots \ar[r] &GW_0^{[n]}\LRf{{\SA}} \ar[r] & K_0({\SA)} \ar[r] & GW_0^{[n+1]}\LRf{{\SA}} \ar[r] &\\
&{\BG}{W}_{-1}^{[n]}\LRf{{\SA}} \ar[r] & 0 \ar[r] & {\BG}{W}_{-1}^{[n+1]}\LRf{{\SA}} \ar[r] &\\
&{\BG}{W}_{-2}^{[n]}\LRf{{\SA}} \ar[r] & 0 \ar[r] & {\BG}{W}_{-2}^{[n+1]}\LRf{{\SA}} \ar[r] &\cdots\\
\enddiagram 
$$
Also consider  the  exact sequence (e.g. \cite[pp 8]{S10})
$$
\diagram
K_0\LRf{{\SA}} \ar[r] & GW_0^{[n+1]} \LRf{{\SA}}  \ar[r] & W^{n+1}\LRf{{\SA}}\ar[r] & 0\\
\enddiagram 
$$
By comparing these two exact sequences, and from the vanishing, it follows: 
$$
\left\{\begin{array}{ll}
{\BG}{W}^{[n]}_{-1}\LRf{{\SA}}=W^{n+1}\LRf{{\SA}} & \forall n\in {\BZ}\\
{\BG}{W}^{[n]}_{-(k+1)}\LRf{{\SA}}\cong {\BG}{W}^{[n+1]}_{-k}\LRf{{\SA}}&\forall n\in {\BZ}, ~k\geq 1\\
\end{array}\right. 
$$ 
So, we have 
$$
\left\{\begin{array}{l}
{\BG}{W}^{[n]}_{-2}\LRf{{\SA}}\cong {\BG}{W}^{[n+1]}_{-1}\LRf{{\SA}}=W^{n+2}\LRf{{\SA}}\\
{\BG}{W}^{[n]}_{-3}\LRf{{\SA}}\cong {\BG}{W}^{[n+1]}_{-2}\LRf{{\SA}}=W^{n+3}\LRf{{\SA}}\\
{\BG}{W}^{[n]}_{-4}\LRf{{\SA}}\cong {\BG}{W}^{[n+1]}_{-3}\LRf{{\SA}}=W^{n+4}\LRf{{\SA}}=W^n\LRf{{\SA}}\\
\end{array}\right. 
$$ 
So, the assertion is valid $\forall n\in {\BZ}$ and $k=1, 2, 3, 4$. 
Now we use induction to complete the proof. Assume the assertion is valid for some $k\geq 1$. Then 
$$
{\BG}{W}^{[n]}_{-(k+1)}\LRf{{\SA}}\cong {\BG}{W}^{[n+1]}_{-(k+1)+1}\LRf{{\SA}}\cong {\BG}{W}^{[n+1]}_{-k}\LRf{{\SA}}\cong W^{n+1+k}\LRf{{\SA}}
$$
\pic $\eop$


We define Karoubi ${\BG}W$ bi spectrum of an exact category with duality, as follows.
\bD{\rm
Let ${\SE}:=\LRf{{\SE}, ^{\vee}, \varpi} \in {\bf eCat}D$ bean exact category with duality.
Assume $1/2\in {\SE}$. On the dg category ${\bf dg}{\SE}$, let $\q$ denote the quasi isomorphisms and $^{\vee}, \varpi$ continue to denote the duality. Define 
$$
\left\{\begin{array}{l}
{\BG}W^{[n]}\LRf{{\SE}}:= {\BG}W^{[n]}\LRf{{\SE}, ^{\vee}, \varpi}:={\BG}W^{[n]}\LRf{{\bf dg}{\SE}, \q, ^{\vee}, \varpi}\in \BiSp\\
{\BG}W^{[n]}_k\LRf{{\SE}}:= {\BG}W^{[n]}_k\LRf{{\SE}, ^{\vee}, \varpi}:
=\pi_k\LRf{\BG}W^{[n]}\LRf{{\SE}}\\
\end{array}\right.
\qquad \forall n, k\in {\BZ}. 
$$
Note
$$
\left\{\begin{array}{ll}
{\BG}W^{[n]}\LRf{{\SE}, ^{\vee}, \varpi}\cong 
{\BG}W^{[n+4]}\LRf{{\SE}, ^{\vee}, \varpi}& \forall n\in {\BZ}\\
{\BG}W^{[2]}\LRf{{\SE}, ^{\vee}, \varpi}\cong 
{\BG}W^{[0]}\LRf{{\SE}, ^{\vee}, -\varpi}& \\
\end{array}\right.
$$
}
\eD

Now we D\'{e}vissage Karoubi  ${\BG}W$ bi spectra.   
\bT[D\'{e}vissage ${\BG}W$ bispectra]\label{BGWDeviss}{\rm 
Let $\LRf{{\SA}, ^{\vee}, \varpi}$ be an  abelian  category with duality, and 
${\SB}\subseteq {\SA}$ be a full subcategory, satisfying the conditions of the setup (\ref{SetupDev}).
Assume that  ${\SA}$ is noetherian, and  every ${\SA}$-module $M$ has a finite ${\SB}$-filtration.
Consider the  inclusion functor: 
$\f: \LRf{{\SB}, ^{\vee}, \varpi} \lra \LRf{{\SA}, ^{\vee}, \varpi}$. 
 Then $\forall~n\in {\BZ}$, the induced map
 $$
  {\BG}{W}^{[n]}\LRf{{\SB}, ^{\vee},  \varpi} \lra  {\BG}W^{[n]}\LRf{{\SA}, ^{\vee},  \varpi}
 $$
 in $\BiSp$ is a homotopy equivalence in $\BiSp$.  This means
 $$
 {\BG}{W}^{[n]}_k\LRf{{\SB}, ^{\vee},  \varpi} \lra  {\BG}{W}^{[n]}_k\LRf{{\SA}, ^{\vee},  \varpi}
 \quad {\rm are~isomorphisms}\quad \forall n, k\in {\BZ}.  
 $$
}
\eT
\pf First, we prove the isomorphisms for the negative groups.
Assume that $k\geq 1$. By (\ref{BGWandWitt335}), along with, \cite[Cor 3.13,  pp 36]{S17}, we have 
$$
{\BG}{W}^{[n]}_{-k}\LRf{{\bf dg}{\SA}, \q, ^{\vee},   \varpi}=W^{n+k}\LRf{{\bf dg}{\SA},  \q,^{\vee}, \varpi}=W^{n+k}\LRf{{\bf D}^b{\SA}, \q, ^{\vee},  \varpi}
$$
So, 
$$
{\BG}{W}^{[n]}_{-k}\LRf{{\bf dg}{\SA}, \q, ^{\vee},  \varpi}=\left\{\begin{array}{lll}
0 & if~n+k=2r+1~odd& by ~\cite[Prop ~5.2]{BW02}\\
W\LRf{{\SA}, ^{\vee}, \varpi} & if~n+k=0 ~mod~4n&by ~ \cite[4.7 ~Thm ~pp~ 376]{B01}\\
W\LRf{{\SA}, ^{\vee}, -\varpi} & if~n+k=2 ~mod~4n&by~\cite[pp ~378]{B01}\\
\end{array}\right.
$$
Likewise
$$
{\BG}{W}^{[n]}_{-k}\LRf{{\bf dg}{\SB}, \q, ^{\vee},   \varpi}=\left\{\begin{array}{ll}
0 & if~n+k=2r+1~odd\\ 
W\LRf{{\SB}, ^{\vee}, \varpi} & if~n+k=0 ~mod~4n\\ 
W\LRf{{\SB}, ^{\vee}, -\varpi} & if~n+k=2 ~mod~4n\\ 
\end{array}\right.
$$
Note that the hypothesis (${\bf B}_0$) of \cite{QSS79} for the inclusion ${\SB}\subseteq {\SA}$ was established in Lemma \ref{Jordan}.
 Comparing these two sets of equations,
 it follows from the D\'{e}vissage theorem for Witt groups   (\ref{HermDev}) 
 that 
$$
{\BG}W^{[n]}_{-k}\LRf{{\bf dg}{\SA}, \q, ^{\vee},  \varpi}
\cong {\BG}W^{[n]}_{-k}\LRf{{\bf dg}{\SA},  \q,  ^{\vee}, \varpi}\qquad \forall n, k\in {\BZ}, ~ k\geq 1. 
$$

Now, we prove the assertion for non negative groups. Let $k\geq 0$. Again, note that the associated triangulated categories 
$\T({\bf dg}{\SA})={\bf D}^b\LRf{{\SA}}$ and $\T({\bf dg}{\SB})={\bf D}^b\LRf{{\SB}}$ are idempotent complete \cite[2.8 Thm]{BS01}. 
Consider the commutative diagram
$$
\diagram
GW_k^{[n]}\LRf{{\bf dg}{\SB}, \q, ^{\vee},  \varpi} \ar[rr] \ar[d] &&{\BG}W^{[n]}_k\LRf{{\bf dg}{\SB}, \q, ^{\vee},  \varpi}\ar[d]\\
GW_k^{[n]}\LRf{{\bf dg}{\SA}, \q, ^{\vee},  \varpi} \ar[rr] && {\BG}W^{[n]}_k\LRf{{\bf dg}{\SA}, \q, ^{\vee},  \varpi}\\
\enddiagram
$$
Then the horizontal arrows are isomorphisms by \cite[Prop 8.7]{S17}, \cite[pp 466]{M23}. We prove that the left vertical arrow is also an isomorphism. 
Note that $Z_0{\bf dg}{\SA}=Ch^b\LRf{{\SA}}$ is an exact category with weak equivalence and duality. 
\bE
\item Let $n=0$. By \cite[pp 1781]{S17}, \cite[pp 450]{M23}
$$
\left\{\begin{array}{l}
GW_k^{[0]}\LRf{{\bf dg}{\SA}, \q,  ^{\vee}, \varpi} =GW_k^{[0]}\LRf{Ch^b{\SA}, \q,  ^{\vee},   +\varpi} \cong GW\LRf{{\SA}, ^{\vee}, \varpi}\\
GW_k^{[0]}\LRf{{\bf dg}{\SB}, \q, ^{\vee},  \varpi} =GW_k^{[0]}\LRf{Ch^b{\SB}, \q,  ^{\vee},   +\varpi}\cong GW\LRf{{\SB}, ^{\vee}, \varpi}\\
\end{array}\right.
$$
By D\'{e}vissage Theorem   \ref{HermDev} of $GW$-groups, we have 
$$
GW\LRf{{\SB}, ^{\vee}, \varpi} \iso
GW\LRf{{\SA}, ^{\vee}, \varpi}
\quad {\rm is~isomorphism}.
$$
So, the first vertical map is isomorphism. 
\item 
Let $n=2$. We have 
$$
\left\{\begin{array}{l}
GW_k^{[2]}\LRf{{\bf dg}{\SA}, \q, ^{\vee},  \varpi} \cong GW_k^{[0]}\LRf{Ch^b{\SA}, \q,  ^{\vee},   -\varpi} \cong GW\LRf{{\SA}, ^{\vee}, -\varpi}\\
GW_k^{[2]}\LRf{{\bf dg}{\SB},  \q, ^{\vee}, \varpi} \cong GW_k^{[0]}\LRf{Ch^b{\SB}, \q, ^{\vee},    -\varpi} \cong GW\LRf{{\SB}, ^{\vee}, -\varpi}\\
\end{array}\right.
$$
Here the first isomorphisms is induced by the shift functor, and latter isomorphisms follows from  \cite[pp 1781]{S17}, \cite[pp 450]{M23}.
By the 
D\'{e}vissage Theorem   \ref{HermDev} of $GW$-groups, we have 
$$
GW\LRf{{\SB}, ^{\vee}, -\varpi} \iso
GW\LRf{{\SA}, ^{\vee}, -\varpi}
\quad {\rm is~isomorphism}.
$$
So, the first vertical map is an isomorphism.
\eE 
 Therefore it follows that the second vertical arrow is also an isomorphism $\forall k\geq 0$ and $n=0, 2 ~mod~4$, as required. 
 
 Now we assume that $n=2m-1$ is odd. Consider the exact triangle \cite[Thm 8.11 pp 1801]{S17} 
 $$
\diagram 
{\BG}W^{[n]}\LRf{\SA} \ar[r] & {\BK}\LRf{\SA} \ar[r] & {\BG}W^{[2m]}\LRf{\SA} \ar[r] & 
S^1\wedge {\BG}W^{[n]}\LRf{\SA}\\
\enddiagram
$$
and the corresponding triangle for ${\SB}$. Comparing the long exact sequences for these two triangles, we obtain the following commutative diagram
of groups:
$$
\diagram 
 K_1({\SB}) \ar[r] \ar[d]^{\wr}& {\BG}W_1^{2m}({\SB}) \ar[r]  \ar[d]^{\wr}
 & {\BG}W^{n}_0({\SB}) \ar[r]  \ar@{-->}[d]^{\wr} & K_0({\SB})  \ar[d]^{\wr}
\ar[r] &  {\BG}W_0^{2m}({\SB}) \ar[r]  \ar[d]^{\wr}& {\BG}W^n_{-1}({\SB}) \ar[d]^{\wr}\\
 K_1({\SA}) \ar[r] & {\BG}W_1^{2m}({\SA}) \ar[r] 
 & {\BG}W^{n}_0({\SA}) \ar[r] & K_0({\SA}) 
\ar[r] &  {\BG}W_0^{2m}({\SA}) \ar[r] & {\BG}W^n_{-1}({\SA})\\
\enddiagram
$$
We have already proved that the solid vertical arrows are isomorphisms. By Five Lemma, the broken vertical arrow is also an isomorphism. For $k\geq 1$ we have the commutative diagram
$$
\diagram 
 K_k({\SB}) \ar[r] \ar[d]^{\wr}& {\BG}W_k^{2m}({\SB}) \ar[r]  \ar[d]^{\wr}
 & {\BG}W^{n}_k({\SB}) \ar[r]  \ar@{-->}[d]^{\wr} & K_{k-1}({\SB})  \ar[d]^{\wr}
\ar[r] &  {\BG}W_{k-1}^{2m}({\SB})  \ar[d]^{\wr}\\
 K_k({\SA}) \ar[r] & {\BG}W_k^{2m}({\SA}) \ar[r] 
 & {\BG}W^{n}_{k-1}({\SA}) \ar[r] & K_{k-1}({\SA}) 
\ar[r] &  {\BG}W_{k-1}^{2m}({\SA})\\ 
\enddiagram
$$
By Five Lemma the broken vertical arrow is an isomorphism. 
So, the proof is complete $\forall n, k\in {\BZ}$. \pic  $\eop$

\section{Application} 
There are two well known applications of the D\'{e}vissage theorem of Quillen 
\cite[pp 107]{M23}, on the ${K}$-theory space ${\bfK}({\SE})$ of small exact categories. 
First, it applies to the category ${\SE}=Coh(X)$ of coherent modules over a noetherian scheme. It follows ${\bfK}\LRf{Coh(X)}\cong
{\bfK}\LRf{Coh(X_{red})}$, where $X_{red}$ denotes the  reduced ssheme of $X$ 
(See \cite[pp 109, 225]{M23}).
In the context of  hermitian theory, this one is not of much interest to us, because 
$Coh(X)$ does not have a natural duality structure. The one of interest to us is the category ${\SE}:=C{\BM}^{V(\m)}\LRf{\spec{R}}$ of  $R$-modules finite length, 
where $(R, \m, \kappa)$ is a regular local ring. As an application of the 
D\'{e}vissage 
theorem it follows from the comments in the introduction (\ref{introDeviss}) that there is a natural homotopy equivalence
 $$
 {\BK}\LRf{C{\BM}^{V(\m)}\LRf{\spec{R}}}\cong {\BK}\LRf{{\SV}(\kappa)}
 $$ 
where ${\SV}(\kappa)$ denotes the category of finite dimensional vector spaces. 
 (see \cite[Cor 3.6.7, pp 111]{M23}). We prove the Hermitian analogue of the same.
 \bT\label{DivCMmXpcmi}{\rm 
 Let $\LRf{R, \m, \kappa}$ be a noetherian regular local ring, $X=\spec{R}$, $Z=V(\m)$, $d=\dim R$. Let 
 $C{\BM}^{Z}\LRf{X}$ denote the category of  $R$-modules $M$ with finite length (note $d=\PDV(M)$). Then 
 $C{\BM}^{Z}\LRf{X}$  has a natural duality $M^{\vee}=Ext^d\LRf{M, R}$. 
 Note $\SV(\kappa)$ also has a natural duality, defined by $M^{\vee}=Hom\LRf{M, \kappa}$. Let $\varpi$ denote the double dual identification, for both. Then there is a natural equivalence, of bi spectra:
 $$
   {\BG}W^{[n]}\LRf{{\SV}(\kappa), ^{\vee},  \varpi}  \iso {\BG}W^{[n]}\LRf{C{\BM}^Z(X), ^{\vee},  \varpi} \qquad in ~\BiSp, \qquad \forall n\in {\BZ}.
 $$
 Consequently,
 $$
  {\BG}{W}^{[n]}_k\LRf{C{\BM}^Z(X), ^{\vee},   \varpi} \cong 
  \left\{\begin{array}{ll}
  W\LRf{{\SV}(\kappa), ^{\vee}, \varpi} &  k\leq -1, n-k-d=0~mod~4\\ 
  0 &  k\leq -1, ~n-k-d=1, 2, 3~mod~4\\ 
  G{W}_k\LRf{{\SV}(\kappa), ^{\vee}, \varpi} & n=0~mod ~4,~ k\geq 0\\
  G{W}_k\LRf{{\SV}(\kappa), ^{\vee}, -\varpi} & n=2~mod ~4, k\geq 0\\
    \end{array}\right.
 $$
 For $k\geq 0$ and $n=2m+1$  odd, we have an exact sequence, as in the proof (\ref{exactKappaOdd}).
 }
 \eT
\pf 
As usual. $\q$ will denote the quasi isomorphisms in the respective categories. 
Let ${\CS}\subseteq C{\BM}^Z(X)$ denote the full  subcategory of semi simple $R$-modules. 
It follows  (\ref{BGWDeviss}) that 
$$
{\BG}W^{[n]}\LRf{{\bf dg}{\CS}, \q, ^{\vee},  \varpi} \lra {\BG}W^{[n]}\LRf{{\bf dg}C{\BM}^Z(X), \q, ^{\vee},  \varpi} 
$$
is a homotopy equivalence. 
Consider the natural functor 
$$
\Phi: \LRf{ {\SV}(\kappa), ^{\vee}, \varpi} \lra \LRf{ {\CS}, ^{\vee}, \varpi}
\quad {\rm where}\quad 
\left\{\begin{array}{ll}
\Phi\LRf{V} = V &  V\in \SV\LRf{\kappa}\\
V^{\vee}=Hom(V, \kappa) & V\in {\SV}(\kappa)\\
M^{\vee}=Ext^d(M, R) & M\in {\CS}\\
\end{array}\right.
$$
Then $\Phi$ induces a duality preserving functor $F:{\bf dg}{\SV}(\kappa)\lra 
\LRf{{\bf dg}{\CS}}$ and the induced functor of associated triangulated categories  
$$
\LRf{{\bfD}^b\LRf{{\SV}(\kappa)} ,  \q,  ^{\vee},\varpi}\iso 
\LRf{{\bfD}^b\LRf{{\CS}} ,  \q, ^{\vee}, \varpi}\qquad \qquad \forall n\in {\BZ}.
$$
is an equivalence \cite[Lem 6.3]{BW02}.
Here 
$^{\vee}$, 
$\varpi$ are induced. 
So, $\Phi$ induces a duality preserving equivalence of the associated triangualted categories. By \cite[pp 466]{M23}, the induced map 
$$
 {\BG}W^{[n]}\LRf{{\bf dg}{\SV}(\kappa), \q, ^{\vee},  \varpi} \lra  {\BG}W^{[n]}\LRf{{\bf dg}C{\BM}^Z(X), \q, ^{\vee},  \varpi} 
 \qquad in \quad \BiSp
 $$
 is an equivalence in $\BiSp$, for all $n\in {\BZ}$. This completes the proof of the main statement. 
 By Prop ~\ref{BGWandWitt335} and remark \ref{notaBalTriW}(\ref{pointFour}), $\forall n, k\in {\BZ}, k \geq 1$ we have 
 $$
  {\BG}{W}^{[n]}_{-k}\LRf{{\SV}(\kappa), ^{\vee}, \varpi} =W^{n+k}\LRf{{\SV}(\kappa), ^{\vee}, 
  \varpi} 
  =
  \left\{\begin{array}{ll}
  W\LRf{{\SV}(\kappa), ^{\vee} 
  \varpi} & if~n+k-d=0~mod~ 4\\
  0 & otherwise\\
  \end{array}\right.
  $$
 This completes the proof regarding the  negative groups. 
  Now, assume $k\geq 0$. By \cite[pp 1781]{S17}, \cite[pp 450]{M23} 
$$
  {\BG}{W}^{[n]}_k\LRf{{\bf dg}{\SV}(\kappa), ^{\vee}, \q, \varpi} =
    {G}{\CW}^{[n]}_k\LRf{{\bf dg}{\SV}(\kappa), ^{\vee}, \q, \varpi} 
  =\left\{\begin{array}{ll}
  {\bf GW}_k\LRf{{\SV}, ^{\vee}, \varpi}& n=0 ~mod~4 \\
  {\bf GW}_k\LRf{{\SV}, ^{\vee}, -\varpi}& n=2 ~mod~4 \\
  \end{array}\right.
  $$
  For odd shift $n=4m+i$,  with $i=1, 3$, we have a long exact sequence
    \begin{equation}\label{exactKappaOdd}
  \diagram 
    {\bf GW}_k\LRf{{\SV}, ^{\vee}, \delta_i \varpi}\ar[r] & K_{k}\LRf{{\SV}(\kappa)} \ar[r] &  G{\CW}^{[n]}_k\LRf{{\bf dg}{\SV}(\kappa), ^{\vee}, \q, \varpi}
  \ar[r] &  {\bf GW}_{k-1}\LRf{{\SV}, ^{\vee}, \delta_i \varpi}\\
  \enddiagram 
  \end{equation}
  where $\delta_1=1$, $\delta_3=-1$. 
  This follows from the exact triangle \cite[pp 1784]{S17}, which is the $G{\CW}$-analogue  of (\ref{extDelts}).
  \pic $\eop$ 
  %

\begin{remark}{\rm 
In deed, the the structure \ref{DivCMmXpcmi} of ${\BG}W^{[n]}_k(C{\BM}^Z(X))$ was instrumental in computing the structure of ${\BG}W$ theory of quasi projective schemes, punctured by a regular point $\m \in X$, which was computed in \cite{M24}. This was the driving motivation for this article. The following (\ref{summHApcmi}) is a simple version of the same.
}
\end{remark}


\bC\label{summHApcmi}{\rm 
Consider the setup of ( \ref{DivCMmXpcmi}).
 Let $\overline{n-d}= n-d ~mod~4$, and $-3 \leq 
\overline{n-d} \leq 0$. Then
\bE
\item For $\overline{n-d}=0$ we have 

\bE
\item  The left side exact sequence terminates, at degree zero:
$$
\diagram 
\cdots \ar[r] &GW_0^{[n]}\LRf{R} \ar[r] & GW_0^{[n]}(X) \ar[r] & GW_0^{n]}(U) \ar[r] &0 \\
\enddiagram 
$$
\item There are isomorphisms
$$
\left\{\begin{array}{l}
GW_{-1-4p}^{[n]}(X) \iso GW_{-1-4p}^{[n]}(U) \\
 GW_{-2-4p}^{[n]}(X) \iso GW_{-2-4p}^{[n]}(U) \\
\end{array}\right. \qquad \forall p\geq 0.
$$
\item There are  exact sequences
$$
\diagram 
0 \ar[r] & GW_{-3-4p}^{[n]}(X) \ar[r] & GW_{-3-4p}^{[n]}(U)   
\ar[r] &\\
W\LRf{\kappa} \ar[r] & GW_{-4-4p}^{[n]}(X) \ar[r] & GW_{-4-4p}^{[n]}(U)  
\ar[r] &0\\
\enddiagram \qquad \forall p\geq 0
$$
\eE 
\item For $\overline{n-d}=\TCP{-1, -2, -3}$ we 
leave it to the readers to write down   similar statements.
\eE
 
}
\eC

\pf Assume $n-d=0 ~mod~4$. We have an exact sequence
$$
\diagram 
\cdots \ar[r] &GW_0^{[n]}\LRf{C{\BM}^Z(X)} \ar[r] & GW_0^{[n]}(X) \ar[r] & GW_0^{[n]}(U)  \\
\ar[r] &0 \ar[r] & GW_{-1}^{[n]}(X) \ar[r] & GW_{-1}^{[n]}(U)   \\
\ar[r] &0 \ar[r] & GW_{-2}^{[n]}(X) \ar[r] & GW_{-2}^{[n]}(U)   \\
\ar[r] &0 \ar[r] & GW_{-3}^{[n]}(X) \ar[r] & GW_{-3}^{[n]}(U)   \\
\ar[r] &W\LRf{\kappa} \ar[r] & GW_{-4}^{[n]}(X) \ar[r] & GW_{-4}^{[n]}(U)  \\
\ar[r] &0\ar[r]& GW_{-5}^{[n]}(X) \ar[r] & GW_{-5}^{[n]}(U)  \\
\ar[r] &0\ar[r] & GW_{-6}^{[n]}(X) \ar[r] & GW_{-6}^{[n]}(U)  \\
\enddiagram
$$
Rest of the proof follows by chasing this exact sequence. \pic $\eop$



\end{document}